\documentclass[amssymb,amsfonts]{amsart}
\usepackage{amssymb, latexsym,amscd}
\def\be#1{\begin{equation} \label{#1}}
\def\bi{\begin{itemize}}
\def\bs{\begin{split}}
\def\es{\end{split}}
\def\ba{\begin{align}}
\def\bas{\begin{align*}}
\def\ea{\end{align}}
\def\eas{\end{align*}}
\def\v{{\tilde v}}
\def\w{\tilde w}
\def\vapr{\tilde v_{\text{appr}}}
\def\u{{\tilde u}}
\def\Im{{\hbox{Im}}}

\def\Re{{\hbox{Re}}}

\def\C{{\hbox{\bf C}}}
\def\R{{{\Bbb R}}}
\def\Hil{{\Bbb H}}

\def\B{{\hbox{\bf B}}}

\def\M{{\hbox{\bf M}}}

\def\Z{{{\Bbb Z}}}
\def\T{{{\Bbb T}}}

\def\eps{\varepsilon}

\def\emph#1{{\it #1}}
\def\textbf#1{{\bf #1}}

\theoremstyle{plain}
\newtheorem{theorem}{Theorem}
\newtheorem{definition}[theorem]{Definition}
\newtheorem{proposition}[theorem]{Proposition}
\newtheorem{lemma}[theorem]{Lemma}
\newtheorem{condition}[theorem]{Condition}
\newtheorem{corollary}[theorem]{Corollary}

\numberwithin{equation}{section}
\numberwithin{theorem}{section}

\include{psfig}

\begin{document}

\title[Symplectic non-squeezing for the KdV flow]{Symplectic non-squeezing of the KdV flow}
\author{J.~Colliander}
\thanks{J.C. was supported in part by  N.S.E.R.C. Grant RGPIN
  250233-03 and the Sloan Foundation.}
\address{University of Toronto}

\author{M.~Keel}
\thanks{M.K. was supported in part by N.S.F. Grant DMS
                         9801558}
\address{University of Minnesota}
\author{G.~Staffilani}
\thanks{G.S. was supported in part by N.S.F. Grant DMS 0100345 and by a grant
from  the Sloan Foundation.}
\address{M.I.T.}
\author{H.~Takaoka}
\address{Kobe University and the University of Chicago}
\thanks{H.T. was supported in part by J.S.P.S. Grant No. 13740087.}
\author{T.~Tao}
\thanks{T.T. was a Clay Prize Fellow and was supported in part by grants
from the Packard Foundation.}
\address{University of California, Los Angeles}

\vspace{-0.3in}
\begin{abstract}
We prove two finite dimensional approximation
results and a symplectic non-squeezing property for the Korteweg-de Vries (KdV) flow on the circle
$\T$.  The nonsqueezing result relies on the aforementioned approximations and
the finite-dimensional nonsqueezing theorem of Gromov \cite{gromov}.
Unlike the work of Kuksin \cite{kuksin} which initiated the investigation of  non-squeezing results
for infinite dimensional Hamiltonian systems, the nonsqueezing argument here does not
construct a capacity directly.
In this way our results are similar to those obtained for the NLS flow by Bourgain \cite{borg:nonsqueeze}.
A major difficulty here though is the lack of any sort of
smoothing estimate which would allow us to easily approximate the infinite dimensional KdV flow by
a finite-dimensional Hamiltonian flow. To resolve this problem we invert the Miura
transform and work on the level of the modified KdV (mKdV) equation, for which smoothing estimates
can be established.
\end{abstract}
\maketitle

\tableofcontents

\section{Introduction}\label{introduction-sec}

This paper is concerned with the symplectic behavior of the Korteweg-de Vries (KdV) flow \be{kdv}
u_t + u_{xxx} = 6 u u_x; \quad u(0,x) = u_0(x)
\end{equation}
on the circle $x \in \T := \R/2\pi\Z$, where $u(t,x)$ is real-valued. In particular we investigate
how the flows may (or may not) be accurately approximated by certain finite-dimensional models,
and then use such an  approximation to conclude a symplectic non-squeezing property.  In order to
describe the symplectic space involved, and state the result precisely, we need to set notation
and recall some previous results describing the well-posedness of the initial value problem
\eqref{kdv}.

On the circle we have the spatial Fourier transform
\begin{align}   \widehat u(k) & := \frac{1}{2\pi}
\int_0^{2\pi} u(x) \exp(-ikx)\ dx \label{spatialft}
\end{align}
for all $k \in \Z$, and the spatial Sobolev
spaces
$$ \| u \|_{H^s_x} := (2\pi)^{1/2} \| \langle k \rangle^s \widehat u \|_{l^2_k}$$
for $s \in \R$, where $\langle k \rangle := (1 + |k|^2)^{1/2}$.  These are natural spaces
for analyzing the KdV flow.

Let $P_0$ denote the mean operator
$$ P_0 u := \frac{1}{2\pi} \int_0^{2\pi} u$$
or equivalently
$$ \widehat{P_0 u}(k) = \chi_{k=0} \widehat u(k).$$
The KdV flow is mean-preserving, and it will be convenient to work in the case when $u$ has mean
zero\footnote{One can easily pass from the mean zero case to the general mean case by a Galilean
transformation $u(t,x) \rightarrow u(t, x - P_0(u)t) - P_0(u)$.}. Accordingly we define the
mean-zero periodic Sobolev spaces $H^s_0$ by
$$ H^s_0 := \{ u \in H^s_x: P_0 u = 0 \}$$
endowed with the same norm as $H^s_x$.

Recent work on the local and global well-posedness theory in $H^s_0$ for \eqref{kdv} is basic to
our results here.
For example,  the geometric conclusions from
finite-dimensional Hamiltonian dynamics which we ultimately need for our nonsqueezing result  can only be
applied in the
setting of rather rough solutions to the initial value problem \eqref{kdv}. We now pause to
summarize some of the analytical techniques that have been developed for the study of such rough
solutions, and the resulting regularity theory (see e.g.  
\cite{borg:xsb},
\cite{kpv:kdv}, 
\cite{cct}, and \cite{ckstt:2}, \cite{ckstt:gKdV}).

\subsection{Summary of local and global well-posedness theory}

If the initial datum $u_0$ for \eqref{kdv} is smooth, then there is a global smooth
solution\footnote{This result can also be obtained by inverse scattering methods, since the KdV
equation is completely integrable.  However, our methods here do not use inverse scattering
techniques, although the special algebraic structure of KdV (in particular, the Miura transform
\cite{MiuraTransform}) is certainly exploited.} $u(t)$ (see e.g. \cite{early}).
We can thus define the non-linear flow map $S_{KdV}(t)$ on $C^\infty(\T)$ by $S_{KdV}(t)
u_0 := u(t)$. In particular this map is densely defined on every Sobolev space $H^s_0$.

If $s \geq -1/2$, then the equation \eqref{kdv} is globally well-posed in $H^s_0$. In other words,
the flow map $S_{KdV}(t)$ is uniformly continuous (indeed, it is analytic) on $H^s_0$ for times
$t$ restricted to a compact interval $[-T,T]$, and for such $s$ we have bounds of the form
\be{kdv-bound} \sup_{|t| \leq T} \| S_{KdV}(t) u_0 \|_{H^s_0} \leq C(s, T, \| u_0 \|_{H^s_0}),
\end{equation}
(see \cite{kpv:kdv}, \cite{ckstt:2}, \cite{ckstt:gKdV} (and also Section \ref{kdv-bound-sec}
below)). For $s < -1/2$ the flow map $S_{KdV}(t)$ is no longer uniformly continuous \cite{cct}
(see also \cite{kpv:counter}) 
or analytic \cite{borg:measures}, so from the
point of view which requires a uniformly continous flow in time, the Sobolev space
$H^{-1/2}_0$ is the endpoint space for the KdV flow. 
Coincidentally, this space is also a
natural phase space for which KdV becomes a Hamiltonian flow; we will have more to say about this
at the end of the introduction. Note however that if one  asks only that the flow be continuous in
time, then global well-posedness for \eqref{kdv} has been established for all $s \geq -1$
in \cite{kappeler} using inverse scattering methods. 
Combining mapping properties of the Miura Transform and the result in \cite{TakaokaTsutsumi}, 
local well-posedness of \eqref{kdv} in $H^s_0$ with a (not uniformly) continuous flow map holds for $-5/8 < s < -1/2$. 

To obtain many of the local and global well-posedness results mentioned above, one iterates
in a certain spacetime Banach space $Y^s$ (defined in \eqref{y-def} below; this space
is a variant of the $X^{s,b}$ spaces used for instance in
\cite{borg:xsb}, \cite{kpv:kdv}), 
which
has the same regularity as $H^s$ in the sense that one has the embedding\footnote{In this paper we
use $A \lesssim B$ to denote an estimate of the form $A \leq CB$, where the implicit constant $C$
may depend on certain parameters such as $s$ which we will specify later in the paper.  Similarly,
$A \ll B$ denotes $B \geq C A$ for some such universal constant $C$.}
$$ \| u \|_{L^\infty_t H^s_x} \lesssim \| u \|_{Y^s}.$$
The nonlinearity is then placed in a companion space $Z^s$ (see \eqref{z-def} below), which is
related to $Y^s$ via an energy estimate of the form
$$ \| \eta(t) u \|_{Y^s} \lesssim \| u(t_0) \|_{H^s} + \| u_t + u_{xxx} \|_{Z^s},$$
for any time $t_0$, and any bump function  $\eta$ supported near $t_0$. (We will elaborate more
upon these spaces and estimates in Section \ref{notation-sec}). The local well-posedness
theory\footnote{Strictly speaking, in order to handle large initial data one must also generalize
this estimate to circles $\R/2\pi\lambda \Z$ of arbitrarily large period, in order to apply
rescaling arguments to make the data small again.  See \cite{ckstt:2}, \cite{ckstt:gKdV}, or
Section \ref{kdv-bound-sec}.} for the KdV equation \eqref{kdv} then hinges on the bilinear
estimate \be{bil} \| (uv)_x \|_{Z^s} \lesssim \| u \|_{Y^s} \| v \|_{Y^s}
\end{equation}
whenever $u, v$ are mean-zero functions and $s \geq -1/2$ (see \cite{kpv:kdv}, \cite{ckstt:2},
\cite{ckstt:gKdV}).

To pass from local well-posedness to global well-posedness one needs to obtain long-time
bounds on the $H^s_0$ norm.  For $-\frac{1}{2} \leq s < 0 $, this has been achieved by means of
the ``$I$-method'', constructing an almost conserved quantity comparable to the $H^s$ norm; see
\cite{ckstt:2}, \cite{ckstt:gKdV}, or Section \ref{kdv-bound-sec}.

\subsection{Low frequency approximation of KdV}

The KdV flow \eqref{kdv} is, formally at least, a Hamiltonian flow on an infinite-dimensional
space.  In order to rigorously apply results from symplectic geometry, we must approximate this
infinite-dimensional flow by a finite-dimensional flow.  Furthermore, in order to apply these
geometric tools, we need that the finite-dimensional flow is itself Hamiltonian.

We begin with a negative result.  Suppose that we wish to study the KdV flow for data $u_0$ whose
Fourier transform is supported on $[-N,N]$ for some large fixed $N$, and specifically to
approximate the KdV flow by a finite-dimensional model.  A first guess for such a model might be
the flow 
\be{pn-kdv} u_t + u_{xxx} = P_{\leq N}( 6 u u_x); \quad u(0) = u_0,
\end{equation}
where $P_{\leq N}$ is the Fourier projection to frequencies $\leq N$:
$$ \widehat{P_{\leq N} u}(k) = \chi_{|k| \leq N} \widehat u(k).$$
Denote the flow map associated to \eqref{pn-kdv} by $S_{P_{\leq N} KdV}(t)$.  This flow has
several advantageous properties; for instance, $S_{P_{\leq N} KdV}(t)$ is a symplectomorphism on
the space $P_{\leq N} H^{-1/2}_0$, associated with a 
natural symplectic structure (see next
subsection).  Since $P_{\leq N} H^{-1/2}_0$ is a finite dimensional space, it is easy to see (e.g.
using $L^2$ norm conservation and Picard iteration) that this flow $S_{P_{\leq N} KdV}$ is
globally smooth and well-defined.  In \cite{borg:nonsqueeze}, the NLS flow $iu_t + u_{xx} = |u|^2
u$ was similarly truncated, and it was shown that the truncated flow was a good approximation to
the original (infinite dimensional) flow.  Unfortunately, the same result does not apply for KdV:

\begin{theorem}\label{kdv-lowfreq-uhoh}  Let $k_0 \in \Z^*$, $T > 0$, $A > 0$.  Then for any $N \gg C(A, T, k_0)$
 there exists initial data $u_0$ with $\| u_0 \|_{H^{-1/2}_0} \leq A$ and $\text{supp}\,(\widehat{u_0}) \subset \{ |k| \leq N \}$
 such that
\be{fluctuate} |\widehat{(S_{KdV}(T) u_0)}(k_0) - \widehat{(S_{P_{\leq N} KdV}(T) u_0)}(k_0)| \geq
c(T,A, k_0)
\end{equation}
for some $c(T,A, k_0) > 0$.
\end{theorem}

In other words, $S_{P_{\leq N} KdV}$ does not converge to $S_{KdV}$ even in a weak topology.

We prove this negative result in Section \ref{counter-sec}.  Basically, the problem is that the
multiplier $\chi_{[-N,N]}$ corresponding to $P_{\leq N}$ is very rough, and this creates
significant deviations between $S_{KdV}$ and $S_{P_{\leq N} KdV}$ near the Fourier modes $k = \pm
N$.  In cubic equations such as mKdV (see \eqref{mkdv-drift} below) or the cubic nonlinear
Schr\"odinger equation, these deviations would stay near the high frequencies
$\pm N$, but in the quadratic KdV equation these deviations create significant fluctuations
near the frequency origin, eventually leading to failure of weak convergence in \eqref{fluctuate}.

Of course there are several obvious ways to modify the finite-dimensional flow \eqref{pn-kdv} in
an attempt to find an effective approximation to the KdV flow for data with Fourier transform
supported on $[-N,N]$, but at least a little bit of care is needed when considering these
modifications.
We let $b(k)$ be the restriction to the integers of a real even bump function adapted
to $[-N, N]$ which equals 1 on
$[-N/2, N/2]$, and consider the evolution 
\be{b-pn-kdv} u_t + u_{xxx} = B(6 u u_x); \quad u(0) =
u_0
\end{equation}
where
$$ \widehat{Bu}(k) = b(k) \widehat{u}(k).$$
Let $S_{B KdV}$ denote the flow map associated to \eqref{b-pn-kdv}.
Observe that this is a finite-dimensional flow on the space $P_{\leq N} H^s_0$.  Unfortunately,
$S_{B KdV}$ is not a symplectomorphism, but we will
explain in \eqref{hiccupcure} below how by conjugating a flow of the form \eqref{b-pn-kdv}
with a simple multiplier operator we will arrive at our desired finite dimensional
symplectomorphism on $P_{\leq N} H^{-\frac{1}{2}}(\T)$
that well-approximates the full KdV flow at low frequencies.  This desired
symplectomorphism is labelled $S^{(N)}_{KdV}(t)$ in \eqref{hiccupcure} below\footnote{The equation which 
defines this flow is given in \eqref{blub} below.}, and once
the aforementioned approximation properties are established, the nonsqueezing result
will follow almost immediately after quoting the finite dimensional nonsqueezing
result of Gromov \cite{gromov}.

The first step in the argument is to show we can approximate
$S_{KdV}$ by $S_{B KdV}$ in the strong $H^s_0$ topology:

\begin{theorem}\label{kdv-approx}
Fix $s \geq -1/2$, $T > 0$, and $N \gg 1$.  Let $u_0 \in H^s_0$ have Fourier transform supported
in the range $|k| \leq N$.   Then
$$
\sup_{|t| \leq T} \| P_{\leq N^{1/2}} (S_{B KdV} u_0(t) - S_{KdV}(t) u_0) \|_{H^{s}_0} \leq
N^{-\sigma} C(s, T, \| u_0 \|_{H^s_0})$$ for some $\sigma = \sigma(s) > 0$.
\end{theorem}

In particular, we can accurately model the KdV evolution for band-limited initial data by a
finite-dimensional flow, at least for frequencies $|k| \leq N^{1/2}$.

The well-posedness statement $\eqref{kdv-bound}$ gives Theorem \ref{kdv-approx} for
all $0 \leq N \leq C(s,T,\| u_0 \|_{H^s_0}, \| \tilde u_0 \|_{H^s_0})$, hence our proof needs only
to consider  $N \geq C(s,T,\| u_0 \|_{H^s_0}, \| \tilde u_0 \|_{H^s_0})$.  This turns out to be the most
interesting case from the point of view of the nonsqueezing applications  of this approximation theorem
which we take up below.

Theorem \ref{kdv-approx} can be viewed as a statement that one can (smoothly) truncate the KdV
evolution at the high frequencies without causing serious disruption to the low frequencies, in
spite of the obstruction posed by Theorem \ref{kdv-lowfreq-uhoh}.  Our second main result (proven in
Section \ref{low-kdv-sec}) is in a similar vein:

\begin{theorem}\label{low-perturb-kdv} Fix $s \geq -1/2$, $T > 0, N \geq 1$.  Let $u_0,
\tilde u_0 \in H^s_0$ be such that $P_{\leq 2N} u_0 = P_{\leq 2N} \tilde u_0$ (i.e. $u_0$ and $\tilde
u_0$ agree at low frequencies).  Then we have,
$$
\sup_{|t| \leq T} \| P_{\leq N}(S_{KdV}(t) \tilde u_0 - S_{KdV}(t) u_0) \|_{H^{s}_0} \leq
N^{-\sigma} C(s, T, \| u_0 \|_{H^s_0}, \| \tilde u_0 \|_{H^s_0})$$ for some $\sigma = \sigma(s) >
0$.
\end{theorem}

By the same reasoning made following Theorem \ref{kdv-approx}, we may
assume in the proof of Theorem \ref{low-perturb-kdv} that
 $N \geq C(s,T,\| u_0 \|_{H^s_0}, \| \tilde u_0 \|_{H^s_0})$.

The point of Theorem \ref{low-perturb-kdv} is that changes to the initial data at frequencies $\geq 2N$ do
not significantly affect
the solution at frequencies $\leq N$, as measured in the strong $H^s_0$ topology.  This is in
stark contrast to the negative result in Theorem \ref{kdv-lowfreq-uhoh}.  The point is that there
is some delicate cancellative structure in the KdV equation which permits the decoupling of high
and low frequencies, and this structure is destroyed by projecting the KdV equation crudely using
\eqref{pn-kdv}.

To prove Theorem \ref{kdv-approx} and Theorem \ref{low-perturb-kdv}, we shall need to exploit the
subtle cancellation mentioned in the previous paragraph in order to avoid the obstructions arising
from Theorem \ref{kdv-lowfreq-uhoh}. We do not know how to do this working directly with the KdV
flow.
Rather, we are able to prove estimates which explicitly account for
this subtle structure in KdV by using the \emph{Miura transform} $u = \M v$, defined by \be{miura-def} u = \M v
:= v_x + v^2 - P_0(v^2).
\end{equation}
As discovered 
in \cite{MiuraTransform}, this transform allows us to conjugate the KdV flow to the
\emph{modified Korteweg-de Vries} (mKdV) flow \be{mkdv-drift} v_t + v_{xxx} = F(v); \quad v(x,0) =
v_0(x)
\end{equation}
where the non-linearity $F(v)$ is given by \be{F-def} F(v) := 6 (v^2 - P_0(v^2)) v_x.
\end{equation}
The modified KdV equation has slightly better smoothing
    properties{\footnote{See Section  \ref{improvedtrilinearsection}, in particular
    Theorem \ref{fn0-improv-thm}.}} than the ordinary KdV equation,
and 
in addition the process of inverting 
the Miura transform adds one degree of regularity
(from $H^{-1/2}_0$ to $H^{1/2}_0$).  In particular, the types of counterexamples arising in Theorem
\ref{kdv-lowfreq-uhoh} do not appear in the mKdV setting, and by proving a slightly more refined 
trillinear estimate than those found in e.g. \cite{ckstt:gKdV} (see in
    particular Theorem \ref{fn0-improv-thm} below)
we are able to prove the above two
theorems by passing to the mKdV setting using the Miura transform.  Of course, in order to close the argument we will need some efficient
estimates on the invertibility of the Miura transform; we set up these estimates (which may be of
independent interest) in Section \ref{miura-sec}.

\subsection{Application to symplectic non-squeezing}
\label{sec:nosqueezestatement}

We can apply the above approximation results to study the symplectic behavior of KdV in
a 
natural phase space $H^{-1/2}_0(\T)$.  Before doing so, we  recall some context and results
from previous works.  We are following here especially the exposition from \cite{hz_book,
kuksin_book}.

\begin{definition} Consider a pair $(\Hil, \omega)$ where $\omega$ is a symplectic form\footnote{That is, a
nondegenerate, antisymmetric form $\omega: \Hil \times \Hil \rightarrow \C$.  We identify in the usual
way $\Hil$ and it's tangent space $T_x \Hil$ for each $x \in \Hil$.} on the Hilbert space $\Hil$.
We say $(\Hil, \omega)$ is the {\emph{symplectic phase space}} of a PDE with Hamiltonian
$H[u(t)]$ if the PDE can be written in the form,
\begin{align}
\label{hamiltonianflow}
\dot{u}(t) & = J \nabla H[u(t)].
\end{align}
\end{definition}
Here $J$ is an almost complex structure\footnote{That is, a bounded, anti-selfadjoint operator
with $J^2 = -(\text{identity})$.}  on $\Hil$, which is compatible with the Hilbert space inner product $\langle \cdot ,
\cdot \rangle$.  That is, for all $u,v \in \Hil$,
\begin{align}
\label{compatible}
\omega(u,v) & = \langle J u, v \rangle.
\end{align}
The notation $\nabla$ in \eqref{hamiltonianflow} denotes the usual gradient with respect
to the Hilbert space inner product,
\begin{align} \langle v, \nabla H[u] \rangle & \equiv d H[u] (v)  \label{gradient}\\
& \equiv {\frac{d}{d \epsilon}}{\big|_{\epsilon = 0}} H[u + \epsilon v].
\end{align}
One easily checks that an equivalent way to write the PDE corresponding to the Hamiltonian
$H[u(t)]$ in $(\Hil, \omega)$ is
\begin{align}
\label{secondd1}
\dot{u}(t) & = \nabla_\omega H[u(t)]
\end{align}
where the symplectic gradient $\nabla_\omega H[u]$ is defined in analogy with \eqref{gradient},
\begin{align}
\label{secondd2}
\omega(v, \nabla_\omega H[u]) & = dH[u](v).
\end{align}
For example, on the Hilbert space $H^{-\frac{1}{2}}_0(\T)$, we can define the symplectic form
\be{symplectic}
\omega_{-\frac{1}{2}}(u,v) := \int_\T u(x) \partial_x^{-1} v(x)\ dx
\end{equation}
where $\partial_x^{-1}: H^{-1/2}_0(\T) \to H^{1/2}_0(\T)$ is the
inverse to 
the differential operator $\partial_x$ defined via the Fourier transform by
$$ \widehat{ \partial_x^{-1} f}(k) := \frac{1}{ik} \widehat{f}(k).$$
The KdV flow \eqref{kdv} is then \emph{formally} the Hamiltonian equation in $(H^{-1/2}_0(\T), \omega_{-\frac{1}{2}})$
corresponding to the (densely
defined) Hamiltonian
\be{hamil-def} H[u] := \int_\T \frac{1}{2} u_x^2 + u^3 dx.
\end{equation}
Indeed, working formally\footnote{By the word `formally', we  mean here that no attempt is made to justify
various differentiations or integration by parts.  Later, when we localize the space
$H^{-\frac{1}{2}}_0$ and Hamiltonian in frequency and write down the corresponding equations, the reader
can carry out the analogous computation where the justification of the necessary calculus will be evident.} we have
for any $v \in H^{-\frac{1}{2}}_0(\T)$,  
\begin{align*}
{\frac{d}{d \epsilon}}{\big|_{\epsilon = 0}} H[u + \epsilon v] & = \int_{\T} u_x v_x + 3 u^2v dx \\
& = \int_\T (-u_{xx} + 3u^2) v dx \\
& = \int_\T \partial_x^{-1} (-u_{xxx} + 6 u u_x) v dx \\
& = -\int_{\T} (-u_{xxx} + 6 u u_x) \partial_x^{-1} v dx \\
& = \omega_{-\frac{1}{2}}(u_{xxx} - 6 u u_x, v) \\
& = \omega_{-\frac{1}{2}}(v, -u_{xxx} + 6 u u_x).
\end{align*}
Comparing \eqref{secondd1}-\eqref{secondd2} with \eqref{kdv}, we see KdV is indeed the Hamiltonian
PDE corresponding to $H[u]$ on the infinite dimensional symplectic space
$(H^{-\frac{1}{2}}_0, \omega_{-\frac{1}{2}})$.  In particular, the flow maps
$S_{KdV}(t)$ are, formally, symplectomorphisms on $H^{-1/2}_0(\T)$.

That the KdV flow arises as a Hamiltonian flow from a symplectic structure as described above 
was discovered by Gardner and Zakharov-Faddeev (see  \cite{gardner, faddeev}).   
A second structure was given by Magri \cite{magri} using $\int u^2 dx$ as Hamiltonian, but it is not 
as convenient as the first structure for our strategy to prove nonsqueezing.  
Roughly speaking, it seems the symplectic 
form in this second structure could possibly be used to establish a nonsqueezing 
property - in the $H^{-\frac{3}{2}}$ topology - of 
a finite dimensional analog of \eqref{kdv}.  However, since the well-posedness theory, and the accompanying estimates,
for the full KdV flow do not presently exist at such rough norms,  we do not see how we could approximate
the full KdV flow in a space as rough as $H^{-\frac{3}{2}}$ with a finite dimensional flow.  The first structure described 
above allows us to adopt this strategy in the space $H^{-\frac{1}{2}}_0$, within which we do 
have well-posedness.  (See below for references for this approach to proving nonsqueezing for PDE.  See 
e.g \cite{olver, dickey} for more details and history of the various symplectic structures for KdV.)

For any $u_* \in H^{-1/2}_0(\T)$, $r > 0$, $k_0 \in \Z^*$, and $z \in \C$, we consider the
infinite-dimensional ball
$$ \B^\infty(u_*; r) := \{ u \in H^{-1/2}_0(\T): \| u - u_* \|_{H^{-1/2}_0} \leq r \}$$
and the infinite-dimensional cylinder
$$ \C^\infty_{k_0}(z; r) := \{ u \in H^{-1/2}_0(\T): |k_0|^{-1/2} |\widehat u(k_0) - z| \leq r \}.$$

The final result of this paper is the following symplectic
non-squeezing theorem,

\begin{theorem}\label{main}  Let $0 < r < R$, $u_* \in H^{-1/2}_0(\T)$, $k_0 \in \Z^*$, $z \in \C$,
and $T>0$.  Then
$$ S_{KdV} (T)(\B^\infty(u_*; R)) \not \subseteq \C^\infty_{k_0}(z; r).$$
In other words, there exists a global $H^{-1/2}_0(\T)$ solution $u$ to \eqref{kdv} such that
$$ \| u(0) - u_* \|_{H^{-1/2}_0} \leq R$$
and
$$ |k_0|^{-1/2} |\widehat{u(T)}(k_0) - z| > r.$$
\end{theorem}

Note that no smallness conditions are imposed on $u_*$, $R$, $z$, or $T$.

Roughly speaking, this Theorem asserts that the KdV flow cannot squash a large ball into a thin
cylinder.  Notice that the balls and cylinders can be arbitrarily far away from the origin, and
the time $T$ can also be arbitrary.  Note though that this result is interesting  even for $u_* = 0, z=0$
and smooth initial data $u_0$, as it tells us that the flow cannot at any time uniformly squeeze the
ball $B^\infty (0, R)$ even at a fixed frequency $k_0$.
By Theorem \eqref{main}, the well-posedness theory
for KdV reviewed above, and density considerations,we know that for any $T, r < R$, there will be some initial
data $u_0 \in  B^{\infty}(0;R)$  for which\footnote{We are using here the statement of the Theorem only
in the case $u_* = 0, z = 0$.  Of course one gets a similar conclusion to the one we draw here, but with different
weights and a different initial data set, by simply using the $L^2$ conservation and time reversability properties
of the flow.  That is, for any $R > r$, there is data $\tilde{u}_0 \in \{\|f\|_{L^2(\T)} \, \leq \, R \} $
such that the evolution $\tilde{u}$ of this data satisfies $|\widehat{\tilde{u}}(k_0,T)| > r.$}
$|\widehat{u}(k_0,T)| > |k_0|^{\frac{1}{2}} r$.  (See \cite{borg:book},
page 96 for the same discussion in the context of a nonlinear Klein-Gordon equation.)
A second immediate application of Theorem \ref{main} to smooth solutions
was highlighted in a different context already in \cite{kuksin}, namely that such smooth
solutions of \eqref{kdv}  cannot uniformly approach
some asymptotic state:  for any neighborhood $B^{\infty}(u_0;R)$ of
the initial data in $H^{-\frac{1}{2}}(\T)$ and for any time $t$,
the diameter of the set $S_{KdV} (t)(B^{\infty}(u_0;R))$ cannot be less than $R$.

The motivation for Theorem \ref{main}, and an important component of its proof, is the finite-dimensional
nonsqueezing theorem of Gromov \cite{gromov} (see also subsequent extensions
in \cite{hz}, \cite{hz_book}).  The extension to the
infinite-dimensional setting provided by a nonlinear PDE seems nontrivial.  The program was
initiated by Kuksin \cite{kuksin}, \cite{kuksin_book} for certain equations where the nonlinear flow is a
compact perturbation of the linear flow.   
That the KdV equation doesn't meet this requirement can be seen by an argument
involving simple computations similar to those supporting Theorem \ref{kdv-lowfreq-uhoh} which
are detailed in
Section \ref{counter-sec} below:  fix $\sigma \ll 1$ and for each integer $N \geq 1$ consider
initial data,
\begin{align*}
u_{0,N}(x) & := \sigma N^{\frac{1}{2}} \cos(Nx).
\end{align*}
Clearly the set $\{u_{0,N} : N = 1,2, \ldots \}$ is bounded in $H^{-\frac{1}{2}}_0$.  However,
when one computes the second iterate\footnote{See in particular equation \eqref{blahblah} for the
notation used here, and if necessary Section \ref{counter-sec} for what we hope is a
sufficiently detailed discussion to allow the reader to reproduce the elementary computations we
quote here.} $u^{[2]}_N$ one sees that it differs from the linear
evolution of $\widehat{u_N^{[0]}}$ at frequency $k = N$ in that,
\begin{align}
\label{onegeorge}
\widehat{u^{[2]}_N}(N,t) - \widehat{u^{[0]}_N}(N,t) & \sim N^{\frac{1}{2}} \sigma^3 e^{i N^3 t}.
\end{align}
By the local well-posedness theory we know, assuming $\sigma$ is sufficiently small compared to
$t$, that the difference between the second iterate and the actual nonlinear evolution $u_N(t)$
of the data $u_{0,N}$ satisfies,
\begin{align}
\label{twogeorge}
\|u_N(t) - u^{[2]}_N(t) \|_{H^{-\frac{1}{2}}_0(\T)} & \lesssim \sigma^4.
\end{align}
Together, \eqref{onegeorge} and \eqref{twogeorge} show that if $\{N_k\}$ is a sequence of integers
relatively prime to one-another\footnote{Note  (for example by examining the iterates and using well-posedness)
that $\widehat{u}_N(t)$ is supported only at frequencies which are integer multiples of $N$.}, then
\begin{align*}
\widehat{u_{N_k}}(N_l,t) - \widehat{u^{[0]}_{N_k}}(N_l,t) & \sim \delta_{k,l}\cdot \sigma^3 \cdot
N_k^{\frac{1}{2}} e^{i N_k^3 t}.
\end{align*}
Hence the set $\{u_{N_k}(t) - u^{[0]}_{N_k}(t) \}$ has no limit point in $H^{-\frac{1}{2}}_0(\T)$.

The  nonsqueezing results of Kuksin were extended to certain
stronger nonlinearities by Bourgain \cite{borg:nonsqueeze,borg:book} - for instance \cite{borg:nonsqueeze}
treats the the cubic non-linear Schr\"odinger  flow on $L^2(\T)$.  In these works, the full solution
map is shown to be well-approximated by a finite dimensional flow constructed by  cutting the solution off to
frequencies $|k| \leq N$ for
some large $N$.   The nonsqueezing results in \cite{borg:nonsqueeze, borg:book} follow then
from a direct application of Gromov's finite dimensional
nonsqueezing result to this approximate flow.

The argument we follow here for the KdV flow is similar
to the work in \cite{borg:nonsqueeze, borg:book}, but seems to require a bit more care.  The complication
seems to us to be somehow rooted in the counterexample of Theorem \ref{kdv-lowfreq-uhoh}, which
clearly exhibits that a sharp cut-off is not appropriate in constructing the approximating flow, but which
seems also to be subtly related to the fact that the estimates necessary to approximate the full KdV flow by a more
gradually truncated flow are unavailable to us when we work directly with the KdV equation.
We have already sketched how we will deal with this difficulty (that is, by passing to the modified
KdV equation) in the discussion which followed Theorem \ref{low-perturb-kdv} above.

We now provide some details of the previous paragraph's sketch, in particular we indicate the difficulties that
arise when one tries to repeat the argument in \cite{borg:nonsqueeze,borg:book}.

Let
$N \geq 1$ be an integer. By simply restricting the form $\omega_{-\frac{1}{2}}$, the space
$(P_{\leq N} H^{-1/2}_0(\T), \omega_{-\frac{1}{2}})$ is a $2N$-dimensional real
symplectic space and hence by general arguments (see e.g. Proposition 1
in \cite{hz_book}) is symplectomorphic to the standard space $(\R^{2N}, \omega_0)$.  We will make
explicit use of such an equivalence below:  any $u \in P_{\leq N} H^{-1/2}_0(\T)$ is determined completely by
\begin{equation}
\begin{split}
(\Re (\widehat{u}(1)), \ldots, \Re (\widehat{u}(N)), \Im (\widehat{u}(1)), \ldots, \Im (\widehat{u}(N))) \\ \equiv
(e_1(u), \ldots, e_n(u), f_1(u), \ldots , f_N(u)) \in \R^{2N}.\end{split} \label{coordinates}
\end{equation}
In terms of the coordinates \eqref{coordinates} the form $\omega_{-\frac{1}{2}}$ defined in
\eqref{symplectic} can be written using the Plancherel theorem as,
\begin{align*}
\omega_{-\frac{1}{2}}(u,v) & = \sum_{
\begin{array}{c}
k = -N \\
k \neq 0
\end{array}}^N
\widehat{u}(-k) \frac{1}{ik} \widehat{v}(k) \\
& = \sum_{k=1}^N \frac{1}{ik} \left( \widehat{u}(-k) \widehat{v}(k) - \widehat{u}(k) \widehat{v}(-k) \right) \\
& = \sum_{k = 1}^N \frac{2}{k} ( \Im (\widehat{v}(k) \overline{ \widehat{u}(k)} )) \\
& = \sum_{k=1}^N \frac{2}{k} ( e_k(u) \cdot f_k(v) - e_k(v) \cdot f_k(u)).
\end{align*}
Write $\Gamma$ for the $N \times N$ matrix $\Gamma \equiv \text{diag}(1, \frac{1}{\sqrt{2}},
\frac{1}{\sqrt{3}}, \ldots \frac{1}{\sqrt{N}})$, $\Lambda \equiv \text{diag}(\Gamma,\Gamma)$, and
$u = (\vec{e}(u), \vec{f}(u)) \in \R^{2N}$
for the coordinates in $P_{\leq N} H^{-1/2}_0(\T)$, we summarize the discussion above by saying,
\begin{align}
\label{equiv_summary}
\omega_{-\frac{1}{2}}(u,v) & = \omega_0(\Lambda(\vec{e}(u), \vec{f}(u)), \Lambda(\vec{e}(v), \vec{f}(v))),
\end{align}
where as before we've written $\omega_0$ for the standard symplectic form on $\R^{2N}$.  In other
words,
\begin{align*}
\Lambda & : (P_{\leq N} H^{-1/2}_0(\T), \omega_{-\frac{1}{2}}) \rightarrow (R^{2N}, \omega_0)
\end{align*}
is a symplectomorphism.

Following \cite{borg:nonsqueeze}, our goal is to find a flow which satisfies three
conditions:  it should be finite dimensional - that is, map $P_{\leq N} H^{-\frac{1}{2}}(\T)$ into itself; it
should be a symplectic map for each $t$; and it should well-approximate the full flow $S_{KdV}(t)$ in a
sense that we will make rigorous momentarily.   For now, we write $S^{(N)}_{\text{Good!}}(t)$ for this
yet to be determined flow.
\begin{equation}
\label{pushover}
\begin{CD}
(P_{\leq N} H^{-\frac{1}{2}}_0, \omega_{-\frac{1}{2}}) @>\Lambda>> (\R^{2N}, \omega_0) \\
@V S^{(N)}_{\text{Good!}}(t)VV     \\
(P_{\leq N} H_0^{-\frac{1}{2}}, \omega_{-\frac{1}{2}}) @>>\Lambda>  (\R^{2N}, \omega_0)
\end{CD}
\end{equation}
Note then that the map,
\begin{align}
\label{longmap}
\Lambda \circ S^{(N)}_{\text{Good!}}(t) \circ \Lambda^{-1} & : (\R^{2N}, \omega_0) \longrightarrow (\R^{2N}, \omega_0)
\end{align}
is likewise a symplectomorphism to which we can apply the finite dimensional theory of
symplectic capacity (see \cite{gromov}, and e.g. \cite{hz_book}). One defines, for any $\vec{x}_* \in \R^{2N},
u^{(N)}_* \in P_{\leq N}
H^{-1/2}_0(\T)$, $r > 0$, $0 < |k_0| \leq N$, and $z \in \C$, the finite-dimensional balls
in $P_{\leq N} H^{-1/2}_0(\T), R^{2N}$, respectively, by the notation,
\begin{align} \label{BN-def}
 \B^N(u^{(N)}_*; r) & := \{ u^{(N)} \in P_{\leq N} H^{-1/2}_0(\T): \| u^{(N)} - u^{(N)}_* \|_{H^{-1/2}_0}
\leq r \} \\
B(\vec{x}_*, r) &:= \{ \vec{x} \in \R^{2N} : |\vec{x} - \vec{x}_*| \leq r \}.
\end{align}
and the finite-dimensional cylinders in the same spaces by,
\begin{align*} \C^N_{k_0}(z; r) &:= \{ u^{(N)} \in P_{\leq N} H^{-1/2}_0(\T): |k_0|^{-1/2} |\widehat{u^{(N)}}(k_0) - z|
 \leq r \} \\
\C_{k_0}(z;r) & := \{(\vec{e}, \vec{f}) \in \R^{2N} :  |(e_{k_0} + \sqrt{-1} f_{k_0}) - z| \leq r) \}.
 \end{align*}

From \cite{gromov}, (see also e.g. Theorem 1, Page 55 in the exposition \cite{hz_book}) we have
the finite-dimensional analogue of Theorem \ref{main}:

\begin{theorem}[\cite{gromov}] \label{gromovtheorem}  Assume that for some $R,r \geq 0, z \in \C, 0 \leq k_0 \leq N,
\vec{x}_* \in \R^{2N}$ there is a symplectomorphism $\phi$ defined
on $B(\vec{x}_*, R) \subset (\R^{2N}, \omega_0)$ so that
$$ \phi(B(\vec{x}_*, R)) \subset C_{k_0}(z;r).$$
Then necessarily $r \geq R$.
\end{theorem}
We apply this theorem to the symplectomorphism $\Lambda \circ S^{(N)}_{\text{Good!}} \circ
\Lambda^{-1}$ defined in \eqref{longmap} above to conclude,

\begin{theorem}\label{finite}
Let $N \geq 1$, $0 < r < R$, $u^{(N)}_* \in P_{\leq N} H^{-1/2}_0(\T)$, $0 < |k_0| \leq N$, $z \in
\C$, and $T > 0$.  Let $S^{(N)}_{\text{Good!}}(T): P_{\leq N} H^{-1/2}_0(\T) \rightarrow P_{\leq N} H^{-1/2}_0(\T)$
be any symplectomorphism.   Then
$$ S^{(N)}_{\text{Good!}}(T)(\B^N(u^{(N)}_*; R)) \not \subseteq \C^N_{k_0}(z; r).$$
\end{theorem}

To deduce Theorem \ref{main} from Theorem \ref{finite}, one would like
to let $N \to \infty$ and show that the flow $S^{(N)}_{\text{Good!}}(T)$
converged to $S_{KdV} (T)$ in some weak sense.  More precisely, one would
need,

\begin{condition}\label{kdv-lowfreq}  Let $k_0 \in \Z^*$, $T > 0 $, $A > 0$, $0 < \eps \ll 1$.
Then there exists an $N_0 = N_0(k_0,T,\eps,A) > |k_0|$ such that
$$ |k_0|^{-1/2} |\widehat{S_{KdV}(T) u_0}(k_0) - \widehat{S^{(N)}_{\text{Good!}}(T) u_0}(k_0)| \ll \eps$$
for all $N \geq N_0$ and all $u_0 \in \B^N(0,A)$.
\end{condition}

Once we find a finite dimensional symplectic flow $S^{(N)}_{\text{Good!}}(t)$ for which Condition \ref{kdv-lowfreq} holds, it is
an easy matter to conclude Theorem
\ref{main}.  Indeed, let $r, R, u_*, k_0, z, T$ be as in that Theorem, and choose $0 < \eps <
(R-r)/2$.  The ball $\B^\infty(u_*; R)$ is contained in some ball $\B^\infty(0; A)$ centered at
the origin.  We choose $N \geq N_0(k_0,T,\eps,A)$ so large that $\| u_* - P_{\leq N} u_*
\|_{H^{-1/2}_0} \leq \eps$.  From Theorem \ref{finite} we can find initial data $u^{(N)}_0 \in P_{\leq N} H^{-\frac{1}{2}}(\T)$
satisfying
$ \| u^{(N)}_0 - P_{\leq N} u_* \|_{H^{-1/2}_0} \leq R - \eps$, and hence by the triangle
inequality,
$$\| u^{(N)}_0 -  u_* \|_{H^{-1/2}_0} \leq R,$$
and so that at time $T$ we have,
$$ |k_0|^{-1/2} |\widehat{S^{(N)}_{\text{Good!}}(T)u^{(N)}_0}(k_0) - z| > r + \eps.$$
If we then apply Condition \ref{kdv-lowfreq} and the triangle inequality we obtain Theorem
\ref{main} with $u_0 := u^{(N)}_0$,
\begin{multline*}
  |k_0|^{-1/2}|z - \widehat{S_{KdV} (T) u^{(N)}_0}(k_0)|  \geq \\
 |k_0|^{-1/2} \left| |z - \widehat{S^{(N)}_{\text{Good!}}(T)u^{(N)}_0}(k_0)|
- |\widehat{S_{KdV} (T) u^{(N)}_0}(k_0) - \widehat{S^{(N)}_{\text{Good!}}(T)u^{(N)}_0}(k_0) | \right|\\
 > r + \epsilon - \epsilon  \quad = r.
\end{multline*}

It remains to define the flow $S^{(N)}_{\text{Good!}}(t)$.  One might first try to follow Bourgain's treatment
of several different
Hamiltonion PDE, notably
the cubic NLS flow on $L^2(\T)$ (see \cite{borg:nonsqueeze}, \cite{borg:book}).   Note that the
Hamiltonian $H[u]$ \eqref{hamil-def} is well-defined on  $(P_{\leq N} H^{-1/2}_0(\T), \omega_{-\frac{1}{2}})$,
and the equation giving the corresponding Hamiltonian flow on this space
can be computed as before to be \eqref{pn-kdv}, which can be viewed either as a PDE or as a system of $2N$
ODE.  The maps $S_{P_{\leq N} KdV}(t)$ are therefore symplectomorphisms, but
from Theorem \ref{kdv-lowfreq-uhoh} we know that Condition \ref{kdv-lowfreq}
fails.

We proceed instead by using a flow of the form \eqref{b-pn-kdv} as follows:
Theorem \ref{kdv-approx} tells us that for any multiplier $\tilde{B}$ of the form
described in \eqref{b-pn-kdv}, the finite dimensional flow $S_{\tilde B KdV}$ provides a good approximation to the
low frequency behavior of $KdV$  However,
the flows $S_{\tilde B KdV}$ are not symplectomorphisms, and hence cannot be candidates for our
flow $S^{(N)}_{\text{Good!}}(t)$ in the discussion above.
Fortunately, there is a quick cure for this hiccup using
the approximation given by Theorem \ref{low-perturb-kdv} as follows:  we will define a symplectic,
finite dimensional flow $S^{(N)}_{KdV}(t)$ on $P_{\leq N}H^{-\frac{1}{2}}_0$ so that the following
diagram commutes.
\begin{equation}
\label{hiccupcure}
\begin{CD}
u_0 \in P_{\leq N} H^{-\frac{1}{2}}_0 @>B>> Bu_0 \\
@V S^{(N)}_{KdV}(t)VV @VV S_{B^2 KdV} (t) V \\
S^{(N)}_{KdV}(t) u_0 @>>B> w(t)
\end{CD}
\end{equation}
We write  explicitly the PDE defining this flow in \eqref{blub} below.
To show that $S^{(N)}_{KdV}(t)$ well approximates $S_{KdV}(t)$ at frequency $k_0$, and hence
qualifies as our choice of $S^{(N)}_{\text{Good!}}(t)$,
we will simply spell out the following:  Theorem \ref{low-perturb-kdv} allows
us to replace $S_{B^2 KdV}(t)$ on the right side of \eqref{hiccupcure} with $S_{KdV}(t)$; and
our choice $N \gg |k_0|$ allows us to ignore both the mappings on the top
of \eqref{hiccupcure} (again, by Theorem \ref{low-perturb-kdv}) and the bottom
of \eqref{hiccupcure} (by the definition of $B$, this is the identity at frequency $k_0$).
We give the details in section \ref{finale-sec} below.

\noindent{\bf{Acknowledgements.}}
This work was conducted at UCLA.  The authors are indebted to Tom Mrowka for his detailed
explanation of symplectic non-squeezing.

\section{Inverting the Miura transform}\label{miura-sec}

As described in the introduction above, our work here on the KdV equation relies on the continuity and invertibility
properties of the Miura transform $u = \M v$, where $\M$ is defined by (see \cite{MiuraTransform}),
$$ \M v := v_x + v^2 - P_0(v^2).$$
The additional $P_0(v^2)$ term here is necessary to make the mean of
$\M v$ vanish.
Let $S_{mKdV}(t)$ denote the flow associated to the mKdV equation \eqref{mkdv-drift}.  Then we
have the intertwining relationship \be{intertwining} \M S_{mKdV}(t) = S_{KdV}(t) \M.
\end{equation}
To see this, we suppose that $v$ solves the mKdV equation \eqref{mkdv-drift}, and set $u := \M v$.  Then
one easily checks,
\bas
u_t + u_{xxx} - 6 u u_x & = (\partial_x + 2v) v_t + (\partial_x + 2v) v_{xxx} + 6 v_x v_{xx} \\
&- 6 (v_x + v^2 - P_0(v^2)) ( v_{xx} + 2vv_x ) \\
& = (\partial_x + 2v) (v_t + v_{xxx} - 6 v^2 v_x + 6 P_0(v^2) v_x ) \\
& = 0.
\end{align*}

Heuristically, the Miura transform acts like a derivative operator $\partial_x$, and in particular
we expect it to be a locally bilipschitz bijection from $H^s_0$ to $H^{s-1}_0$.  The purpose of
this section is to make this heuristic rigorous for the range $s \geq
1/2$.  (See also \cite{kt}, which studies
the Miura transform for the larger range $s \geq 0$.)

In what follows we shall make frequent use of the well-known Sobolev multiplication law
\be{sob-mult} \| uv\|_{H^s(\T)} \lesssim \| u \|_{H^{s_1}(\T)} \| v \|_{H^{s_2}(\T)}
\end{equation}
whenever $s \leq \min(s_1,s_2)$ and $s \leq s_1 + s_2 - \frac{1}{2}$, with at least one of the two
inequalities being strict.

From \eqref{sob-mult} it is clear that $\M$ is a locally Lipschitz\footnote{By this we mean that
$\M$ is Lipschitz on every ball in $H^s_0(\T)$, with a Lipschitz constant depending on the ball.}
map from $H^s_0(\T)$ to $H^{s-1}_0(\T)$ for $s \geq 1/2$ (in fact $s > 0$ would suffice).  The
main result of this section is to invert this statement:

\begin{theorem}\label{miura-invert}  Let $s \geq 1/2$.  Then the map $\M$ is a bijection from $H^s_0(\T)$
to $H^{s-1}_0(\T)$, and the inverse map $\M^{-1}$ is a locally Lipschitz map from $H^{s-1}_0(\T)$
to $H^s_0(\T)$.
\end{theorem}

\begin{proof}
We shall focus on the endpoint case $s=1/2$.  We shall see at the end of the proof that the higher
regularity cases $s>1/2$ then follow from the endpoint case and standard elliptic regularity
theory.  We remark that the arguments here (based on a variational approach) are unrelated to the
rest of the paper and can be read independently.

Since the linearization $v \mapsto v_x$ of the Miura transform $\M$ is clearly bilipschitz from
$H^{1/2}_0(\T)$ to $H^{-1/2}_0(\T)$ it is tempting to treat the lower order terms $v^2 - P_0(v^2)$
as perturbations to be iterated away.  This works well if $v$ and $\M(v)$ are small, however for
large $v$ it appears that iterative techniques alone cannot obtain this result\footnote{However,
iterative techniques do allow us to bootstrap low regularity estimates to high regularity
estimates, basically because $\M$ is elliptic and $v$ lies above the critical regularity
$H^{-1/2}$ for $\M$ (and for mKdV).  The strategy of this argument will be to use variational
estimates to obtain a preliminary estimate in very rough norms, and use iteration to improve this
to estimates in the correct norms $H^{1/2}_0(\T)$ and $H^{-1/2}_0(\T)$.}.  Indeed, we shall need
to also rely on variational techniques, and in particular we will use the well-known connection
between the Miura transform and the spectral theory of Schr\"odinger operators.  The key identity
here is 
\be{m-connect} (\frac{d}{dx} + v)(-\frac{d}{dx} + v) = - \frac{d^2}{dx^2} + (v_x + v^2) =
-\frac{d^2}{dx^2} + \M(v) + P_0(v^2).
\end{equation}

We shall work entirely with the smooth functions in $H^{1/2}_0(\T)$ and $H^{-1/2}_0(\T)$, and
obtain bilipschitz bounds for $\M$ on these functions; it will then be clear from standard
limiting arguments that one has bilipschitz bounds in general.

Let $u \in H^{-1/2}_0(\T)$ be smooth.  We consider the problem of finding a smooth function $v \in
H^{1/2}_0(\T)$ with $\M v = u$, showing this $v$ is unique and of estimating $v$ in terms of $u$.
This will be achieved by studying the self-adjoint Schr\"odinger operator $L = L_u$ defined by
$$ L := -\frac{d^2}{dx^2} + u(x)$$
and the associated energy functional $E[\phi] = E_u[\phi]$ defined on $H^1(\T)$ by
$$ E[\phi] := \langle L \phi, \phi \rangle = \int_{\T} \phi_x^2(x) + u(x) \phi^2(x)\ dx.$$
Since $L$ is a self-adjoint elliptic operator on a compact manifold $\T$, it has a discrete
spectrum $\lambda_1 \leq \lambda_2 \leq \ldots$ with $\lambda_n \to +\infty$.  In particular we
have a lowest eigenvalue $\lambda_1 = \lambda_1(u) \in \R$, and a non-zero (real-valued)
eigenfunction $\phi_1$ with $L \phi_1 = \lambda_1 \phi_1$.  \emph{A priori} $\phi_1$ is only in
$H^1(\T)$, but since $u$ is smooth one can use the equation $L \phi_1 = \lambda_1 \phi_1$ to
deduce that $\phi_1$ is also smooth.

Our analysis here shall rely solely on $\lambda_1$.  It is interesting to note that the work in
\cite{borg:measures}, which is at a similar level of scaling to $H^{-1/2}_0(\T)$, uses the entire
spectrum $\lambda_n$ of the operator $L$.

From construction of $E[\phi]$ we observe that \be{E-min} E[\phi] \geq \lambda_1 \int_{\T} \phi^2
\end{equation}
for all $\phi \in H^1(\T)$, with equality attained if and only if $\phi$ is a $\lambda_1$
eigenfunction of $L$. (As we shall see, $\lambda_1$ is an isolated eigenvalue, so equality only
occurs when $\phi = c \phi_1$ for some $c$.)  Thus $\lambda_1$ can be described in a variational
manner.

Since $u \in H^{-1/2}_0(\T)$ we see that $E[1] = 0$, thus $\lambda_1$ must be non-positive.  If $u
\not \equiv 0$ then 1 is not an eigenfunction, and so $\lambda_1$ becomes strictly negative.

We now claim that $\phi_1$ cannot vanish anywhere.  If it had a double zero at some point, i.e.
$\phi_1(x_0) = \partial_x \phi_1(x_0) = 0$, then from the second-order ODE $L \phi_1 = \lambda_1
\phi_1$ and the Picard existence theorem for ODE we see that $\phi_1 \equiv 0$, a contradiction.
Now suppose that $\phi_1$ had a simple zero at $x_0$, so in particular $\phi_1$ changed sign.  Let
$\phi_1 = \phi_1^+ + \phi_1^-$ denote the positive and negative components of $\phi_1$.  An
integration by parts shows that
$$ E[\phi_1^+] = \int_{\phi_1 > 0} L \phi_1(x) \phi_1(x) = \lambda_1 \int \phi_1^+(x)^2.$$
This implies that $\phi_1^+$ is a $\lambda_1$ eigenfunction of $L$, which contradicts the fact
that all such eigenfunctions are smooth\footnote{Alternatively, one can smooth $\phi_1^+$ at the
zeroes of $\phi_1$ to contradict \eqref{E-min}.}.   Thus $\phi_1$ is nowhere vanishing; without
loss of generality we may take $\phi_1$ to be positive and $L^2$-normalized (which uniquely
identifies $\phi_1$).  If we now define $v$ to be the logarithmic derivative of $\phi_1$
$$ v(x) := \frac{\partial_x \phi_1(x)}{\phi_1(x)}$$
then $v$ is smooth and we have
$$ v_x = \frac{\partial_{xx} \phi_1}{\phi_1} - (\frac{\partial_x \phi_1}{\phi_1})^2 = u -
\lambda_1 - v^2$$
(since $L \phi_1 = \lambda_1 \phi_1$) and hence
$$ u = v_x + v^2 + \lambda_1.$$
Taking means of both sides we see that \be{lv} -\lambda_1 = P_0(v^2)
\end{equation}
and hence $u = \M v$.

This shows existence of $v$ such that $u=\M v$.  Observe from \eqref{m-connect} and an
integration by parts that \be{E-ident} E[\phi] = \int (\phi_x - v \phi)^2\ dx - P_0(v^2) \int
\phi^2;
\end{equation}
from this and \eqref{lv} we immediately see that \eqref{E-min} holds (which we already knew), and
that equality occurs if and only if $\phi_x = v \phi$, or in other words if $\phi$ is a constant
multiple of $\exp( \partial_x^{-1} v )$.  In particular this shows that $v$ is unique, for if we
had $\M v = \M \tilde v$ then the above argument yields $\exp(\partial_x^{-1} v)$ is a constant
multiple of $\exp(\partial_x^{-1} \tilde v)$, which implies $v = \tilde v$ if $v, \tilde v$ both
lie in $H^{1/2}_0(\T)$.

We have now shown that $\M$ is smooth, locally Lipschitz, and bijective on smooth functions with
mean zero. To extend this to $H^{1/2}_0(\T)$ and $H^{-1/2}_0(\T)$ we need some \emph{a priori}
estimates on $\M^{-1}$ in these norms.

Let $u \in H^{-1/2}_0(\T)$ and $v \in H^{1/2}_0(\T)$ be smooth
functions such that $u = \M v$.
For this discussion we will allow implicit constants to depend on the $H^{-1/2}_0(\T)$ norm of
$u$.  Write $U := \partial_x^{-1} u$, thus $\| U \|_{H^{1/2}_0(\T)} \lesssim 1$. We observe from
integration by parts, H\"older, Sobolev, and Gagliardo-Nirenberg that \bas
E[\phi] &= \int \phi_x^2 + \int u \phi^2 \\
&= \| \phi \|_{\dot H^1}^2 - 2 \int U \phi \phi_x \\
&\geq \| \phi \|_{\dot H^1}^2 - C \| U \|_{L^4_x} \| \phi \|_{L^4_x} \| \phi_x \|_{L^2_x} \\
&\geq \| \phi \|_{\dot H^1}^2 - C \| U \|_{H^{1/2}_0(\T)} \| \phi \|_{H^{1/2}_x} \| \phi \|_{\dot H^1} \\
&\geq \| \phi \|_{\dot H^1}^2 - C \| \phi \|_{L^2}^{1/2} \| \phi \|_{H^1}^{1/2} \| \phi \|_{\dot H^1}.
\end{align*}
In particular we have the coercivity bound
$$
E[\phi] + C \| \phi \|_{L^2}^2 \gtrsim \| \phi \|_{H^1}^2$$ for all $\phi \in H^1(\T)$.  Applying this
to $\phi=\phi_1$ in particular and recalling the upper bound on
$\lambda_1$ we obtain the eigenvalue bound \be{l1-bound} -C \leq \lambda_1 \leq
0
\end{equation}
and the preliminary eigenfunction bound
$$ \| \phi_1 \|_{H^1} \lesssim 1.$$
From \eqref{sob-mult} and the $H^{-1/2}_0(\T)$ bound on $u$ we thus have
$$ \| u \phi_1 \|_{H^{-1/2}} \lesssim 1$$
which by the eigenfunction equation $L \phi_1 = \lambda_1 \phi_1$ implies the better eigenfunction
bound \be{better} \| \phi_1 \|_{H^{3/2}} \lesssim 1.
\end{equation}

Now we estimate $v$.  From \eqref{lv} and \eqref{l1-bound} we have the preliminary bound
$$ \| v \|_{L^2} \lesssim 1;$$
since $u = \M v$, we thus have
$$ \| v_x - u \|_{L^1} \lesssim 1.$$
Since $L^1$ and $H^{-1/2}_0(\T)$ both embed into $H^{-3/4}$ (for instance) we thus have by Sobolev
that
$$ \| v \|_{L^4} \lesssim \| v \|_{H^{1/4}_0} \lesssim \| v_x \|_{H^{-3/4}} \lesssim 1.$$
Returning once again to the equation $u = \M v $, we thus have
$$ \| v_x - u \|_{L^2} \lesssim 1$$
which then implies \be{v-bound} \| v \|_{H^{1/2}_0(\T)} \lesssim 1.
\end{equation}
In particular we have
$$ \| \partial_x^{-1} v \|_{L^\infty} \lesssim
\| \partial_x^{-1} v \|_{H^{3/2}_0(\T)} \lesssim  \| v \|_{H^{1/2}_0(\T)} \lesssim 1,$$ and thus
$\exp(\partial_x^{-1} v)$ is bounded above and below.  Since $\phi_1$ is a constant multiple of
$\exp(\partial_x^{-1} v)$, we thus see from \eqref{better} that \be{phi-const} |\phi_1(x)| \sim 1
\hbox{ for all } x \in \T.
\end{equation}

We have obtained good bounds for $v = \M^{-1}(u)$ and for the ground state $\phi_1$.  We now
establish that $\M^{-1}$ is Lipschitz for smooth $v$ in a given bounded subset of $H^{1/2}_0$.
From the inverse function theorem and the fact (from \eqref{sob-mult}) that $\M$ is a locally
uniformly $C^2$ map from $H^{1/2}_0$ to $H^{-1/2}_0$, it suffices to show that the derivative map
$\M'(v): H^{1/2}_0 \to H^{-1/2}_0$ is uniformly invertible for $v$ in this set.

A direct computation shows
$$ \M'(v)(w) = (1-P_0) (\partial_x + 2v) w$$
We shall invert this explicitly.

\begin{lemma}\label{mv-invert}  We have
$$ \M'(v)^{-1} = A[\exp(-2 \partial^{-1}_x v)] \partial_x^{-1} A[\exp(2 \partial^{-1}_x v)]$$
where for any positive function $V \in H^{3/2}(\T)$, $A[V]: H^{\pm 1/2}_0(\T) \to H^{\pm
1/2}_0(\T)$ is the operator
$$ A[V](w) := Vw - \frac{V}{P_0(V)} P_0(Vw).$$
\end{lemma}

We recommend that the reader think of $\M'(v)$ and $\M'(v)^{-1}$ as perturbations of $\partial_x$
and $\partial_x^{-1}$ respectively.

\begin{proof}
We have \bas
\M'(v) &= (1-P_0) \exp(-2\partial^{-1}_x v) \partial_x \exp(2 \partial^{-1}_x v)\\
&= (1-P_0) \exp(-2\partial^{-1}_x v) \partial_x (1-P_0) \exp(2 \partial^{-1}_x v).
\end{align*}
Also, observe that $A[V]$ is the inverse of $(1-P_0)V^{-1}$ on $H^{\pm 1/2}_0(\T)$.  The claim
follows.
\end{proof}

Since $H^{3/2}$ is a Banach algebra (by \eqref{sob-mult}), we have \be{vb} \| \exp(\pm 2
\partial^{-1}_x v) \|_{H^{3/2}} \lesssim \exp( C \| \partial^{-1}_x v \|_{H^{3/2}} ) \lesssim
\exp( C \| v\|_{H^{1/2}_0(\T)} ) \lesssim 1.
\end{equation}
Thus from Lemma \ref{mv-invert} we see that $\M'(v)^{-1}$ is uniformly bounded from $H^{-1/2}_0$
to $H^{1/2}_0$.

We have now proven Theorem \ref{miura-invert} at the endpoint $s=1/2$.  We now sketch how one can
use elliptic regularity theory to bootstrap this to higher
regularities $s>1/2$.

Let us first show the boundedness of $\M^{-1}$ from $H^{s-1}_0$ to $H^s_0$ for smooth functions.
In other words, if $u = \M v$ is smooth, we wish to show that
$$ \| v \|_{H^s_0} \lesssim C( \| u \|_{H^{s-1}_0} ).$$
From the $H^{1/2}$ theory we already know that
$$ \| v \|_{H^{1/2}_0} \lesssim C( \| u \|_{H^{s-1}_0} ).$$
Suppose for the moment that $1/2 < s < 3/2$.  We write \bas
\| v \|_{H^s_0} &\lesssim \| v_x \|_{H^{s-1}_0} \\
&\lesssim \| \M v \|_{H^{s-1}_0} + \| (1-P_0) v^2 \|_{H^{s-1}_0} \\
&\lesssim \| u \|_{H^{s-1}_0} + \| v^2 \|_{H^{s-1}}.
\end{align*}
If $s < 3/2$, then by \eqref{sob-mult} we see that $\| v^2 \|_{H^{s-1}} \lesssim \| v
\|_{H^{1/2}_0}^2 \lesssim C( \| u \|_{H^{s-1}_0})$, which establishes boundedness.  By iterating
this type of argument again one can cover the case $3/2 \leq s < 5/2$, and so forth until we
obtain boundedness for all $s>1/2$.  The local Lipschitz property for $\M^{-1}$ is proven
similarly and is left to the reader.
\end{proof}

From the above Theorem, the analyticity of $\M$, and the inverse function theorem we see in fact
that $\M^{-1}$ is locally uniformly $C^m$ as a map from $H^{s-1}_0(\T)$ to $H^{s}_0(\T)$, for any
integer $m$ and any $s \geq 1/2$.

\section{The Fourier restriction spaces $Y^s$ and $Z^s$}\label{notation-sec}

In view of the results of the last section, we see that to analyze the KdV flow in the $H^{s-1}_0$
topology it will suffice to analyze the mKdV flow in the $H^s_0$ topology.  We now review the
basic machinery (from \cite{borg:xsb}, \cite{kpv:kdv}, \cite{ckstt:2}, \cite{ckstt:gKdV}) for
doing so.

If $u(x,t)$ is a function on the cylinder $\T \times \R$ with mean zero at every time,
and $s,b \in \R$, we
define the $X^{s,b} = X^{s,b}(\T \times \R)$ norm by
$$
\| u \|_{X^{s,b}} := \| \widehat u(k, \tau) \langle k \rangle^s \langle \tau - k^3\rangle^b
\|_{L^2_{\tau,k}}$$ where $L^2_{\tau,k}$ is with respect to Lebesgue measure $d\tau$ in the $\tau$
variable and counting measure in the $k$ variable, $\langle x \rangle^2 \equiv (1 + |x|^2)$, and
the space-time Fourier transform $\widehat
u(k,\tau)$ is given for $k \in \Z^*$, $\tau \in \R$ by
$$ \widehat u(k, \tau) := \int_{\T \times \R} e^{-2 \pi i (x k + t\tau)} u(x,t)\ dx dt.$$
We use the same notation here as for the purely spatial Fourier transform \eqref{spatialft},
relying on context to distinguish the two.

We also need the spaces \be{y-def} \| u \|_{Y^s} := \| u \|_{X^{s,1/2}} + \| \langle k \rangle^s
\widehat u \|_{L^2_k L^1_\tau}
\end{equation}
and \be{z-def} \| u \|_{Z^s} := \| u \|_{X^{s,-1/2}} + \| \frac{\langle k \rangle^s \widehat
u}{\langle \tau - k^3  \rangle} \|_{L^2_k L^1_\tau}.
\end{equation}

Observe that we have the crude estimate \be{crude} \| u \|_{Z^s} \lesssim \| u \|_{X^{s,0}} = \| u
\|_{L^2_t H^s_x}
\end{equation}
which will be useful for controlling quartic or higher order error terms; often we will be
localized in time and just estimate $L^2_t H^s$ by $L^\infty_t H^s$.  Here and in the sequel, we
always allow implicit constants to depend on the exponent $s$.

We can restrict the space $Y^s$ to a time interval $I \subseteq \R$ in the usual manner as
$$ \| u \|_{Y^s_I} := \inf \{ \| v \|_{Y^s}: v|_{\T \times I} = u \}.$$
Similarly we can restrict the $Z^s$ norm.  In practice we shall work in a fixed time interval
(usually $[-T,T]$) and implicitly restrict all of our norms to this interval.

Now we give some embeddings for the $Y^s$ and $Z^s$ spaces.  Since the Fourier transform of an
$L^1$ function is continuous and bounded, we have from \eqref{y-def} that \be{sup-est} Y^s
\subseteq C_t H^s_x \subseteq L^\infty_t H^s_x.
\end{equation}

We have the ``energy estimate'',
\be{energy-est}
\| \eta(t) v \|_{Y^s} \lesssim \| v(t_0) \|_{H^s_0} + \| v_t + v_{xxx} \|_{Z^s}
\end{equation}
for any $t_0 \in \R$ and any bump function $\eta$ supported on $[t_0-C, t_0+C]$.
(\cite{borg:xsb}, see also Lemma 3.1, \cite{ckstt:gKdV}; see
Lemmas 3.1 - 3.3 in \cite{kpv:negative} for analogous estimates in
the nonperiodic context.)

Recall too the main estimate from \cite{ckstt:gKdV} (see Proposition 1 in that paper), namely,
 \be{multi} \|(1 - P_0) \left(((1-P_0) \prod_{i=1}^k u_i) w_x\right) \|_{Z^s} \lesssim (\prod_{i=1}^k \| u \|_{Y^s})
\| w \|_{Y^s}
\end{equation}
for any $s \geq 1/2$ and any integer $k \geq 2$, where the implicit constant depends on $k$.  (We
shall only use \eqref{multi} with $k=2,3,4$).
This particular estimate is crucial (especially at the endpoint $s=1/2$) in order to prove
the local (and global) well-posedness of the modified KdV equation \eqref{mkdv-drift} in
$H^s_0(\T)$ for $s \geq 1/2$.

It would be very convenient if the $Z^s$ on the left-hand side of \eqref{multi} could be replaced by $Z^{s+\sigma}$
for some $\sigma > 0$; this \emph{extra smoothing} estimate would make it easy to ignore the
high-frequency components of the evolution and concentrate on the low frequency evolution.
Unfortunately it is easy to see (by modifying the examples in \cite{kpv:kdv}) that such estimates
fail, especially at $s = 1/2$.  Fortunately, as we will see in the next section there are some other
ways to improve the trilinear version  of \eqref{multi} which will be useful for our approximation results.

\section{An improved trilinear estimate}
\label{improvedtrilinearsection}

The estimate \eqref{multi} with $k = 2$ allows us to estimate the cubic nonlinearity $F(v)$ defined in
\eqref{F-def}.  However for our analysis we shall need a refined version of this estimate.

The first step is to decompose $F$ into ``resonant'' and ``non-resonant'' components. In the
following analysis we shall always assume that $v$ has mean zero.

We start with the Fourier inversion formula
$$ v(x) = \sum_{k \in \Z^*} \widehat v(k) \exp(ikx)$$
for $v \in H^s_0$, where $\Z^* := \Z \backslash \{0\}$ is the set of the non-zero integers. A
direct computation gives that the Fourier transform of $F(v)$ is \be{wfk} \widehat{F(v)}(k) = 6
\sum_{k_1,k_2,k_3 \in \Z^*: k_1+k_2+k_3 = k; k_1+k_2 \neq 0}  \widehat v(k_1) \widehat v(k_2) ik_3 \widehat
v(k_3)
\end{equation}
for all $k \in \Z^*$.  The constraint $k_1 + k_2 \neq 0$ arises since we have subtracted the mean
$P_0(v^2)$ from $v^2$ in the definition of $F(v)$.  Observe
that $F(v)$ is a perfect derivative and so has mean zero and thus no Fourier component at 0.
\begin{lemma}\label{form-lemma}  We have
$$F(v) = F_0(v,v,v) + F_{\neq 0}(v,v,v)$$
where the ``resonant'' trilinear operator $F_0$ is given by \be{f0-def} \widehat{F_0(u,v,w)}(k) :=
-6ik \widehat u(k) \widehat v(k) \widehat w(-k)
\end{equation}
for $k \in \Z^*$, and the ``non-resonant'' trilinear operator $F_{\neq 0}$ is
defined by
\be{fn0-def} \widehat{F_{\neq 0}(u,v,w)}(k) := \; \; - \mspace{-120mu}
\sum_{
\begin{array}{c}
k_1,k_2,k_3 \in \Z^*:\\
 k_1+k_2+k_3 = k;\\
 (k_1+k_2)(k_1+k_3)(k_2+k_3) \neq 0
\end{array}} \mspace{-100mu}
2i (k_1+k_2+k_3) \widehat u(k_1) \widehat v(k_2) \widehat w(k_3)
\end{equation}
for $k \in \Z^*$. \end{lemma}

\begin{proof}
Consider the right-hand side of \eqref{wfk}, and break the sum into pieces according to how many
of the quantities $k_1 + k_3, k_2 + k_3$ are zero.
There is a single term in the sum for which  $k_2+k_3 = k_1+k_3 = 0$, and the summation in this case is $F_0(v,v,v)$.
If just $k_2+k_3$ is zero, then the
total contribution of this case vanishes since the summand in this case is antisymmetric with respect to
swapping $k_2$ and $k_3$.  Similarly if just $k_1+k_3$ is zero.  The remaining portion of the
summation can be seen to be $F_{\neq 0}(v,v,v)$ by a symmetrization in $k_1, k_2, k_3$.
\end{proof}

If $k = k_1+k_2+k_3$, then we have the fundamental \emph{resonance identity} \be{resonance} k^3 -
(k_1^3 + k_2^3 + k_3^3) = 3(k_1+k_2)(k_1+k_3)(k_2+k_3)
\end{equation}
(see e.g. \cite{borg:xsb}).  This justifies the terminology that $F_0$ is ``resonant'' but
$F_{\neq 0}$ is ``non-resonant''.

We remark that, if $u$, $v$, $w$ are real, then $F_0(u,v,w)$ and $F_{\neq 0}(u,v,w)$ are also
real, despite the presence of the imaginary $i$ in the definitions of these quantities.  This
follows from identities such as $\widehat u(-k) = \overline{\widehat u(k)}$. We leave the details to the
reader.  We also remark that eventually these two functions will be estimated in absolute value,
so the constants which appear (e.g. the minus signs out front) will play no role.

\subsection{The $F_0$ (resonant) estimate}

We now give an estimate for $F_0$.   Morally at least, the bound we give follows from the trilinear
version of \eqref{multi},
but we present an independent proof here for the sake of completeness.

\begin{lemma} For any $s \geq 1/2$, and any $u,v,w \in Y^s$ with mean zero, we have
\be{f0} \| F_0(u,v,w) \|_{Z^s} \lesssim \| u \|_{Y^s} \| v \|_{Y^s} \| w \|_{Y^s}.
\end{equation}
\end{lemma}

\begin{proof}
We shall just prove the endpoint case $s=1/2$, as the general case easily follows (e.g. by using
the identity $\partial_x^{s-1/2} F_0(u,v,w) = F_0(\partial_x^{s-1/2}u,v,w)$).

Split $u = \sum_{k \in \Z^*} u_k$, where $u_k$ is a complex-valued function whose spatial Fourier
transform is supported on a single frequency $k$.  Observe that
$$ F_0(u,v,w) = \sum_k F_0(u_k, v_{-k}, w_k).$$
Thus if we show that
$$
\| F_0(u_k,v_{-k},w_k) \|_{Z^{1/2}} \lesssim \| u_k \|_{1/2,1/2} \| v_{-k} \|_{1/2,1/2} \| w_k
\|_{1/2,1/2},$$ then the claim \eqref{f0} follows by summing in $k$ and using Cauchy-Schwartz in
$u$ and $v$ (just estimating the $w_k$ term crudely by $w$).

Fix $k$, and define the function $G_{u_k}(t)$ by  $u_k(x,t) = e^{ikx}
e^{ik^3 t} G_{u_k}(t)$,  
so that $$\|u_k\|_{X^{s, \delta}} = <k>^s \|<\tau>^\delta \widehat{G}_{u_k}(\tau)\|_{L^2_\tau(\R)},$$ similarly for 
$G_{v_{-k}}$ and $G_{w_k}$.
The claim then collapses (after some translation in frequency space) to the one-dimensional
temporal estimate
$$ \|G_{u_k} G_{v_{-k}} G_{w_k} \|_{H^{-1/2}_t} \lesssim 
\| G_{u_k} \|_{H^{1/2}_t} \| G_{v_{-k}} \|_{H^{1/2}_t} \| G_{w_k} \|_{H^{1/2}_t}$$
and
$$ \| \widehat{G_{u_k} G_{v_{-k}} G_{w_k}} \langle \tau \rangle^{-1} \|_{L^1_\tau} \lesssim \| G_{u_k} \|_{H^{1/2}_t}
 \| G_{v_{-k}}\|_{H^{1/2}_t} \| G_{w_k} \|_{H^{1/2}_t}.$$
But both left-hand sides can be estimated by $\| G_{u_k} G_{v_{-k}} G_{w_k} \|_{L^2_t}$, and the claim follows easily
from H\"older and Sobolev.
\end{proof}

\subsection{The $F_{\neq 0}$ (nonresonant) estimate}

We now turn to the non-resonant portion $F_{\neq 0}$ of the non-linearity.  In analogy with
\eqref{multi}, \eqref{f0} we have the estimate \be{fn0} \| F_{\neq 0}(u,v,w) \|_{Z^s} \lesssim \| u
\|_{Y^s} \| v \|_{Y^s} \| w \|_{Y^s}
\end{equation}
for all $s \geq 1/2$ and $u,v,w \in Y^s$ with mean zero.  This estimate can be proven by the
techniques used to prove \eqref{multi} in \cite{ckstt:gKdV}, but we shall obtain it as a consequence
of a slightly stronger version, which we now state.

We first need some Littlewood-Paley notation.  If $N$ is an integer power of two, we let $P_N$
denote the dyadic projection operator
$$ \widehat{P_N u}(k) = \chi_{N \leq |k| < 2N} \widehat u(k).$$
If $N_0,N_1,N_2,N_3$ are four integer powers of two, we let $soprano,alto,tenor,baritone$ be a permutation of the indices
$0,1,2,3$ such that
$$ N_{soprano} \geq N_{alto} \geq N_{tenor} \geq N_{baritone}.$$

\begin{theorem}\label{fn0-improv-thm}  Let $N_0, N_1, N_2, N_3$ be integer powers of two.  Then
\be{fn0-improv} \| P_{N_0} F_{\neq 0}(P_{N_1} u,P_{N_2} v,P_{N_3} w) \|_{Z^{1/2}} \lesssim
(\frac{N_0}{N_{soprano}})^\sigma N_{tenor}^{-\sigma} \| u \|_{Y^{1/2}} \| v \|_{Y^{1/2}} \| w
\|_{Y^{1/2}}
\end{equation}
for some absolute constant\footnote{The quantity $\sigma$ shall vary from line to line.} $\sigma >
0$.
\end{theorem}

This means that \eqref{fn0-improv} is only sharp when the output frequency $N_0$ is essentially
the highest frequencies, and the two lowest frequencies $N_{tenor}$ and $N_{baritone}$ are $O(1)$.
This means that very low Fourier modes can influence high modes, but medium and high modes do not.  In addition, 
the high modes do not have
much influence on the low modes\footnote{
That is, when the soprano and alto dyadic factors were high frequencies and $N_0$ were low, we have a small
first factor on the right side of \eqref{fn0-improv}.}. From \eqref{fn0-improv}
one can easily obtain \eqref{fn0} by summing\footnote{More precisely, one first observes that the
left-hand side of \eqref{fn0-improv} vanishes unless $N_{soprano} \sim N_{alto}$.  Then one
decomposes $u$, $v$, $w$ into dyadic pieces and exploits orthogonality of the projections $P_N$ in
the $Y_s$ and $Z_s$ spaces.  We omit the details.} in the $N_i$.

The estimates \eqref{f0} and \eqref{fn0-improv}
give some intuition for why it's possible to find a finite dimensional approximation to the mKdV flow - and hence,
using the Miura transform, for the KdV flow as well:  the only nonlinear
interactions for which we now have no sharpened estimates are the resonant interactions coming
from $F_0$ (which doesn't mix frequencies) and the high-low-low interactions in $F_{\neq 0}$.
Heuristically, then, we might start believing that if we truncate high frequencies, the evolution
will not see much of a difference at low frequencies. In fact, it is possible to use these estimates
to prove low frequency approximation theorems for mKdV analogous to 
Theorems \ref{kdv-approx}, \ref{low-perturb-kdv}, but we do not write out these results explicitly
in this work.  

The rest of this section is devoted to the proof of Theorem \ref{fn0-improv-thm}.  We remark that the
computations in this section are not needed elsewhere in the paper, and the reader may wish to
take \eqref{fn0-improv} for granted on the first pass and move to the next section.

\begin{proof}  
We begin by reviewing some (non-trivial) estimates from \cite{ckstt:gKdV}.

The proof of \eqref{fn0-improv} relies mainly on the trilinear estimate \be{tri-key} \| u_1 u_2
u_3 \|_{L^2_{x,t}} \lesssim \| u_1 \|_{X^{0,\frac{1}{2} - \frac{1}{100}}} \| u_2
\|_{X^{0,\frac{1}{2}-\frac{1}{100}}} \| u_3 \|_{X^{\frac{1}{2} -
\frac{1}{100},\frac{1}{2}-\frac{1}{100}}}
\end{equation}
proven in Section 7 of \cite{ckstt:gKdV}.  This estimate can be viewed as a trilinear variant of
the $L^6_{x,t}$ Strichartz estimate in \cite{borg:xsb}, and its proof requires a small amount of
elementary number theory.

We will also use the following estimate, which follows relatively quickly from some bounds found in \cite{ckstt:gKdV},
\be{decay} \| \frac{ \langle k \rangle^s \widehat{F_{\neq 0}(u,v,w)}(k)
}{\langle \tau - k^3 \rangle^{1-\delta} } \|_{L^2_k L^1_\tau} \lesssim \| u \|_{Y^s} \| v
\|_{Y^s} \| w \|_{Y^s}
\end{equation}
for all $s \geq 1/2$ and some $\delta > 0$. To establish \eqref{decay}, recall Theorem 3 from \cite{ckstt:gKdV}
\begin{equation}
\label{productestimate} 
\| \prod_{i = 1}^k u_i \|_{X^{s-1, \frac{1}{2}}} \; \lesssim \prod_{i =1}^k \| u_i \|_{Y^s},
\end{equation}
for $s \geq \frac{1}{2}$.  We need also equation (9.2) in \cite{ckstt:gKdV}, which also holds when $s \geq \frac{1}{2}$,
\begin{align}
\label{startitout}
\| \frac{\langle k \rangle^s  \chi_{k \neq 0} ( \widehat{ (1 - P_0)u_1 \cdot (1 - P_0)u_2})(k,\tau)}{\langle \tau - 
k^3 \rangle^{1 - \delta} } \|_{L^2_k L^1_\tau} & \lesssim \|u_1\|_{X^{s - 1, \frac{1}{2}}} \|u_2 \|_{X^{s - 1, \frac{1}{2}}}.
\end{align}
Combining these two and writing for the moment 
\begin{align} 
W(k,\tau) & \equiv \chi_{k \neq 0} (k) ( \widehat{ (1 - P_0)u_x \cdot (1 - P_0)(v w)})(k,\tau)  \nonumber \\
& = \chi_{k \neq 0} (k) \mspace{-30mu}\sum_{\begin{array}{c}
k_1+ k_2 + k_3 = k  \\
k_1, k_2 + k_3 \neq 0
\end{array}}  \mspace{-25mu} i k_1 \widehat{u}(k_1) \widehat{v}(k_2) \widehat{w}(k_3),    \label{Wdef}
\end{align}
we conclude by \eqref{productestimate}
\begin{align}
\| \frac{\langle k \rangle^s W(k,\tau)}{\langle \tau - k^3 \rangle^{1 - \delta}} \|_{L^2_k L^1_\tau} & \lesssim
\|u_x \|_{X^{s - 1, \frac{1}{2}}} \cdot \| vw \|_{X^{s - 1, \frac{1}{2}}} \nonumber \\
& \lesssim \|u\|_{Y^s} \|v\|_{Y^s} \| w \|_{Y^s} \label{middlestep}.
\end{align}
We quickly conclude \eqref{decay} from \eqref{middlestep}:  looking at the definition of the norms involved, 
one sees that without loss of generality we may assume $\widehat{u}(k), \widehat{v}(k), \widehat{w}(k) 
\geq 0$ .  Next, by replacing the factor $\widehat{u}(k)$ appearing in \eqref{middlestep}, \eqref{Wdef} with 
$\chi_{k_1 \geq 0}(k_1) \widehat{u}(k_1), \chi_{k_1 \leq 0} \widehat{u}(k_1)$, one concludes \eqref{middlestep} with the function
$W$ now replaced by, 
\begin{align*}
W_{II}(k,\tau) & \equiv \chi_{k \neq 0} (k) \mspace{-35 mu}\sum_{\begin{array}{c}
k_1+ k_2 + k_3 = k \\
k_1, k_2 + k_3 \neq 0
\end{array}}  \mspace{-40mu} |k_1| \widehat{u}(k_1) \widehat{v}(k_2) \widehat{w}(k_3).
\end{align*}
Repeating this argument while interchanging the roles of $k_1, k_2$, and then $k_1, k_3$ and summing  gives
\eqref{middlestep} with $W$ replaced with,
\begin{align*}
W_{III}(k,\tau) & 
\equiv \chi_{k \neq 0} (k) \mspace{-90mu} \sum_{\begin{array}{c}
k_1 \cdot k_2 \cdot k_3 \neq 0 \\
(k_1+ k_2) (k_2 + k_3) (k_1 + k_3) \neq 0
\end{array}}  \mspace{-90mu} (|k_1| + |k_2| + |k_3|) \widehat{u}(k_1) \widehat{v}(k_2) \widehat{w}(k_3).
\end{align*}
By the definition of $F_{\neq 0}$ \eqref{fn0-def} this yields \eqref{decay}.

We now begin the proof of \eqref{fn0-improv}.  It will suffice to prove the estimate
\be{fn0-main}
\| P_{N_0} F_{\neq 0}(P_{N_1} u,P_{N_2} v,P_{N_3} w) \|_{X^{1/2,-1/2}} \lesssim
(\frac{N_0}{N_{soprano}})^\sigma N_{tenor}^{-\sigma} \| u \|_{X^{1/2,1/2}} \| v \|_{X^{1/2,1/2}}
\| w \|_{X^{1/2,1/2}}.
\end{equation}
Indeed, this estimate already controls the $X^{1/2,-1/2}$ portion of the $Z^{\frac{1}{2}}$ norm.  To control
the $L^2_k L^1_\tau$ portion, we observe  from H\"older that the left-hand side of
\eqref{fn0-main} controls
$$
\| \frac{ \langle k \rangle^{\frac{1}{2}} \widehat{(P_{N_0} F_{\neq 0}(P_{N_1} u,P_{N_2} v,P_{N_3} w)}(k)
}{\langle \tau - k^3 \rangle^{1 + \delta} } \|_{L^2_k L^1_\tau},$$ and the claim follows by a suitable
interpolation with \eqref{decay} (decreasing $\sigma$ if necessary).

It remains to prove \eqref{fn0-main}.  By duality this is equivalent to
$$ |\int\int u_0 \partial_x^{-1} F_{\neq 0}(u_1,u_2,u_3)\ dx dt| \lesssim
(\frac{N_0}{N_{soprano}})^\sigma N_{tenor}^{-\sigma} \| u_0 \|_{X^{-3/2,1/2}} \| u_1
\|_{X^{1/2,1/2}} \| u_2 \|_{X^{1/2,1/2}} \| u_3 \|_{X^{1/2,1/2}}
$$
where $u_i$ has Fourier support on the region $|k_i| \sim N_i$.  We have inserted the
$\partial_x^{-1}$ multiplier to cancel the $(k_1+k_2+k_3)$ factor in \eqref{fn0-def}.

The right-hand side is comparable to
\begin{equation}
\label{stopper}
(\frac{N_0}{N_{soprano}})^\sigma N_{tenor}^{-\sigma} \frac{(N_0 N_1 N_2 N_3)^{1/2}}{N_0^2}
\prod_{j=0}^3 \| u_j \|_{X^{0,1/2}}.
\end{equation}
Note that we may assume $N_{soprano} \sim N_{alto}$ 
since the left-hand side of \eqref{fn0-main}
vanishes otherwise. Hence the right side of \eqref{stopper} is bounded below (throwing away the factor
$\left(\frac{N_0}{N_{soprano}}\right)^\sigma$) by  
$$
N_{tenor}^{1/2-\sigma} N_{baritone}^{1/2} N_{soprano}^{-1} \prod_{j=0}^3 \| u_j
\|_{X^{0,1/2}}
$$
Taking spacetime Fourier transforms and taking advantage of the frequency
localization, we thus reduce to showing 
\be{fn-0targ}
\begin{split}
\mspace{50mu} | \mspace{-80mu}\sum_{\begin{array}{c}
k_0,k_1,k_2,k_3 \in \Z^*\\
 k_0 + k_1+k_2+k_3=0; \\
(k_1+k_2)(k_2+k_3)(k_2+k_3) \neq 0
\end{array}} \mspace{-60mu}\int \prod_{j=0}^3 \widehat{u_j}(k_j,\tau_j) \ d\tau|\\
& \mspace{-30mu}\lesssim N_{tenor}^{1/2-\sigma} N_{baritone}^{1/2} N_{soprano}^{-1} \prod_{j=0}^3 \| u_j
\|_{X^{0,1/2}}
\end{split}
\end{equation}
where $d\tau$ is integration over the three-dimensional space $\{ (\tau_0,\tau_1,\tau_2,\tau_3)
\in \R^4: \tau_0+\tau_1+\tau_2+\tau_3 = 0\}$ with measure $d\tau :=
\delta(\tau_0+\tau_1+\tau_2+\tau_3) \prod_{j=0}^3 d\tau_j$. We remark that the above estimate is
now symmetric with respect to permutations of $k_0,k_1,k_2,k_3$.

Without loss of generality we may assume that the $\widehat{u}_j$ are all non-negative.  The next
step is to exploit the implicit $\langle \tau_j - k_j^3\rangle^{1/2}$ denominators.
  From the fundamental identity \eqref{resonance},
\be{fund} \sum_{j=0,1,2,3} \tau_j - k_j^3 = - \sum_{j=0}^3 k_j^3 = 3 (k_1+k_2)(k_2+k_3)(k_1+k_3)
\end{equation}
we see that \bas
\sup_{j=0,1,2,3} \langle \tau_j - k_j^3 \rangle &\gtrsim |k_1+k_2| |k_2+k_3| |k_1+k_3| \\
&= |k_{soprano}+k_{baritone}| |k_{alto}+k_{baritone}| |k_{tenor}+k_{baritone}|.
\end{align*}

By symmetry we may assume that the supremum on the left-hand side is attained when $j=0$.

\begin{lemma}\label{test-lemma} We have
\be{test} |k_{soprano}+k_{baritone}| |k_{alto}+k_{baritone}| |k_{tenor}+k_{baritone}| N_{baritone}
\gtrsim N_{soprano}^2
\end{equation}
\end{lemma}

\begin{proof}
We have four cases:

{\bf Case 1:}  $N_{baritone} \ll N_{tenor} \ll N_{alto}$.  Then the left-hand side of \eqref{test}
is comparable to $N_{soprano}^2 N_{tenor} N_{baritone}$.

{\bf Case 2:} $N_{baritone} \sim N_{tenor} \ll N_{alto}$.  Then the left-hand side of \eqref{test}
is at least $\gtrsim N_{soprano}^2 N_{tenor}$.

{\bf Case 3:} $N_{baritone} \ll N_{tenor} \sim N_{alto}$.  Then the left-hand side of \eqref{test}
is comparable with $N_{soprano}^3 N_{baritone}$.

{\bf Case 4:} $N_{baritone} \sim N_{tenor} \sim N_{alto}$.  Then at least one of $k_1+k_2$,
$k_2+k_3$, $k_1+k_3$ must have magnitude $\sim N_{soprano}$ (since they sum to $-2k_0$).  Since
the other two factors have magnitude at least 1, the left-hand side of \eqref{test} is $\gtrsim
N_{soprano}^2$.
\end{proof}

From this lemma, we have
$$ \langle \tau_0 - k_0^3 \rangle \gtrsim N_{soprano}^2 N_{baritone}^{-1}.$$
Thus to prove \eqref{fn-0targ} it will suffice to show that
$$
| \sum_{\begin{array}{c}
k_0,k_1,k_2,k_3 \in \Z^*: \\
k_0 + k_1+k_2+k_3=0; \\
(k_1+k_2)(k_2+k_3)(k_2+k_3) \neq 0
\end{array}} \int N_{tenor}^{-1/2+\delta}
\langle \tau_0 - k_0^3 \rangle^{1/2} \prod_{j=0}^3 \widetilde{u_j}(k_j,\tau_j)\ d\tau | \lesssim
\prod_{j=0}^3 \| u_j \|_{X^{0,1/2}}.
$$
At least one of $k_1, k_2, k_3$ is $O(N_{tenor})$; by symmetry let's suppose it's $k_3$.  Then we
can bound $N_{tenor}^{1/2-\delta}$ by $k_3^{1/2-\delta}$, and then by undoing the Fourier
transform and doing some substitutions the estimate becomes
$$ |\int\int v_0 v_1 v_2 v_3\ dx dt| \lesssim
\| v_0 \|_{X^{0,0}} \| v_1 \|_{X^{0,1/2}} \| v_2 \|_{X^{0,1/2}} \| v_3 \|_{X^{1/2-\sigma,1/2}}.$$
But this follows directly from \eqref{tri-key} if $\sigma$ is small enough.  This proves
\eqref{fn0-improv}.
\end{proof}

\section{Proof of Theorem \ref{low-perturb-kdv}:  KdV low frequencies are stable under high
  frequency perturbations of data.}
\label{low-kdv-sec}

We now prove Theorem \ref{low-perturb-kdv}.  Fix $s$, $T$, $u_0$, $\tilde u_0$.

We have no upper bound on the time $T$, and so in particular we cannot hope to control the flow
$S_{KdV}(t)$ on the entire interval $[-T,T]$ by a single application of the local well-posedness
theory.  On the other hand, because of the uniform bounds \eqref{kdv-bound} we see that we can
divide $[-T,T]$ into a bounded number $C(s, T, \| u_0 \|_{H^s_0}, \| \tilde u_0 \|_{H^s_0})$ of
time intervals such that the local well-posedness theory can be used on each interval.  It will
thus suffice to prove a local-in-time version of Theorem \ref{low-perturb-kdv}; more precisely, it
will suffice to show

\begin{proposition}\label{low-perturb-kdv-small} Fix $s \geq -1/2$, and $N' \geq 1$.  Let $u_0,
 \tilde u_0 \in H^s_0$ be such that $P_{\leq N'} u_0 = P_{\leq N'} \tilde u_0$.  Then, if
$T'$ is sufficiently small depending on $s$, $\| u_0 \|_{H^s_0}$, and $\| \tilde u_0 \|_{H^s_0}$,
we have
$$
\sup_{|t| \leq T'} \| P_{\leq N' - (N')^{1/2}}(S_{KdV}(t) \tilde u_0 - S_{KdV}(t) u_0)
\|_{H^{s}_0} \leq (N')^{-\sigma} C(s, \| u_0 \|_{H^s_0}, \| \tilde u_0 \|_{H^s_0})$$ for some
$\sigma = \sigma(s) > 0$.
\end{proposition}

The exponent $1/2$ in $(N')^{1/2}$ is not particularly important here; any exponent between 0 and
1 would suffice.

To see how this proposition implies the theorem, first recall that we may assume that $N$ is
large, $N \geq C(s, T, \| u_0 \|_{H^s_0}, \| \tilde u_0 \|_{H^s_0} )$, since the claim in Theorem
\ref{low-perturb-kdv} trivially follows from \eqref{kdv-bound} otherwise.  (This same remark also
applies of course in Proposition \ref{low-perturb-kdv-small}, allowing us to assume
$N' \geq C(s, \| u_0 \|_{H^s_0}, \| \tilde u_0 \|_{H^s_0})$ there too.)  From \eqref{kdv-bound}
we may divide $[-T,T]$ into $C(s, T, \| u_0 \|_{H^s_0}, \| \tilde u_0 \|_{H^s_0})$ time
intervals, such that on each interval (a time-translated version of) Proposition
\ref{low-perturb-kdv-small} holds.  Consider for example the first such time interval $[0,T']$ on the positive
real axis.  We start with  $N' := 2N$ and apply Proposition \ref{low-perturb-kdv-small}, to get,
\begin{align*}
\sup_{t \in [0,T']} \| P_{\leq N' - (N')^{1/2}}(S_{KdV}(t) \tilde u_0 - S_{KdV}(t) u_0)
\|_{H^{s}_0}& \leq (N')^{-\sigma} C(s, \| u_0 \|_{H^s_0}, \| \tilde u_0 \|_{H^s_0}).
\end{align*}
Before moving on to the next subinterval, modify $S_{KdV}(T') \tilde u_0 $ on frequencies $|k| \leq
N' - (N')^{\frac{1}{2}}$ to agree with $S_{KdV}(T') u_0$.  By the local-well posedness theory
and the triangle inequality, we can proceed as on the first subinterval,  decrementing $N'$ by
$(N')^{1/2}$ each time we apply Proposition \ref{low-perturb-kdv-small}, to obtain Theorem \ref{low-perturb-kdv}
 if $N$ (and hence
$N'$) is sufficiently large.

It remains to prove Proposition \ref{low-perturb-kdv-small}.
Henceforth we allow our implicit constants to depend on $s$,
$\| u_0 \|_{H^s_0}$, and $\| \tilde u_0 \|_{H^s_0}$.

Define,
$$ v_0 := \M^{-1} u_0; \quad v(t) := S_{mKdV}(t) v_0;$$
$$ \tilde v_0 := \M^{-1} \tilde u_0; \quad \tilde v(t) := S_{mKdV}(t) \tilde v_0;$$
from Theorem \ref{miura-invert} we thus have
$$ \| v_0 \|_{H^{s+1}_0}, \| \tilde v_0 \|_{H^{s+1}_0} \leq C$$
while from \eqref{intertwining} we have
$$ S_{KdV}(t) u_0 = \M v(t); \quad S_{KdV}(t) \tilde u_0 = \M \tilde v(t).$$
Our task is thus to show that \be{mv-mv} \sup_{|t| \leq T'} \| P_{\leq N' - (N')^{1/2}}(\M \tilde
v(t) - \M v(t)) \|_{H^{s}_0} \leq C (N')^{-\sigma}.
\end{equation}
Henceforth we allow the quantity $\sigma > 0$ to vary from line to line.

We first investigate the discrepancy between $\tilde v$ and $v$ at time 0.

\begin{lemma}\label{v-disc}  With $v_0, \tilde v_0$ defined as above, we have,
$$ \| P_{\leq N'} (\tilde v_0 - v_0) \|_{H^{s+1}_0} \leq C (N')^{-\sigma}.$$
\end{lemma}

\begin{proof}
From the definitions and our assumptions on $u_0, \tilde u_0$ we have
$$ P_{\leq N'}(\M \tilde v_0 - \M v_0) = 0.$$
On the other hand, from Theorem \ref{miura-invert} we have
$$ \| P_{\leq N'} (\tilde v_0 - v_0) \|_{H^{s+1}_0}
\leq C \| \M P_{\leq N'} \tilde v_0 - \M P_{\leq N'} v_0 \|_{H^s_0}.$$ Thus by the triangle
inequality, it will suffice to show the commutator estimate,
\begin{align} \label{commutator}
 \| \M P_{\leq N'} v_0 - P_{\leq N'} \M v_0 \|_{H^s_0} & \leq C (N')^{-\sigma}
 \end{align}
and similarly for $\tilde v_0$.

Clearly it will suffice just to consider $v_0$.  From the definition \eqref{miura-def} of the tranform
$\M$ and the fact that $P_0$,
$P_{\leq N'}$ and $\partial_x$ all commute, we have \bas \M P_{\leq N'} v_0 - P_{\leq N'} \M v_0
=& (1-P_0) [(P_{\leq N'} v_0)^2 - P_{\leq N'} (v_0)^2] \\
=& (1-P_{\leq N'})[ (P_{\leq N'} v_0)^2 ] \\
&- (P_{\leq N'} - P_0)[ ((1 - P_{\leq N'}) v_0) ((1 + P_{\leq N'}) v_0)].
\end{align*} 
But the last two terms have an $H^s_0$ norm of $O((N')^{-\sigma})$ for some $\sigma > 0$; this can
be seen by the Sobolev multiplication law \eqref{sob-mult}, the $H^{s+1}$ bound on $v_0$, and the estimate
$$ \| (1 - P_{\leq N'}) v \|_{H^s} \lesssim N^{-\sigma} \| v \|_{H^{s+\sigma}}$$
to extract the $(N')^{-\sigma}$ decay from the high frequency projection $1-P_{\leq N'}$.  The
claim follows.
\end{proof}

We still have to prove \eqref{mv-mv}. It will suffice to show that \be{v-v} \sup_{|t| \leq T'} \|
P_{\leq N' - (N')^{1/2}}(\tilde v(t) - v(t)) \|_{H^{s+1}_0} \leq C (N')^{-\sigma}.
\end{equation}
This is basically because the commutator of $\M$ with $P_{\leq N' - (N')^{1/2}}$ is small thanks
to the argument in the proof of Lemma \ref{v-disc}.  We omit the details as they are very similar
to those in Lemma \ref{v-disc}.

From Lemma \ref{v-disc} we see that $\tilde v_0$ and $v_0$ are almost identical at low frequencies
$|k| \leq N'$.  In fact, because the solution map $S_{mKdV}(t)$ is locally
Lipschitz\footnote{Since we are assuming $T'$ to be small this follows directly from the local
well-posedness theory.} in $H^{s+1}_0$, we may assume that
\be{v-proj} P_{\leq N'}(\tilde
v_0 - v_0) = 0,
\end{equation}
since the general case then follows by modifying $\tilde v_0$ (or $v_0$) by a small amount in
$H^{s+1}_0$ and using the Lipschitz property.

Henceforth we assume \eqref{v-proj}, so that the low frequency ($|k| \leq N'$) portions of $\tilde
v(t)$ and $v(t)$ are identical at time 0.  Our task is to prove \eqref{v-v}, which asserts that
the slightly lower frequency ($|k| \leq N' - (N')^{1/2}$) portions of $\tilde v(t)$ and $v(t)$ are
still very close together at later times.  This will be achieved primarily through the improved
trilinear estimate \eqref{fn0-improv}.

In what follows we assume all our spacetime norms are restricted to the time interval $[-T',T']$.

From the local well-posedness theory of mKdV (See\footnote{Strictly
  speaking, when the data  $v_0$, $\tilde v_0$
has large $H^{s+1}_0$ norm one has to first rescale the torus by a suitable scaling parameter
$\lambda$ in order to close the iteration, but this has no significant effect on our argument.
The details are carried out in  \cite{ckstt:gKdV},
  \cite{ckstt:2}.} \eqref{energy-est}, \eqref{f0}, \eqref{fn0}, 
or \cite{borg:xsb}, \cite{kpv:negative}, \cite{ckstt:gKdV})
we have the local
estimates \be{vtv-est} \| v \|_{Y^{s+1}} + \| \tilde v \|_{Y^{s+1}} \leq C
\end{equation}
if the time $T'$ is chosen sufficiently small depending on the $H^{s+1}_0$ norms of $v_0$, $\tilde
v_0$.

The frequency interval $[N'-(N')^{1/2}, N']$ contains $O((N')^{1/4})$ intervals of the form $[M, M
+ (N')^{1/4}]$.  By orthogonality and the pigeonhole principle, we see that there must exist one
of these intervals $[M, M+(N')^{1/4}]$ such that \be{med-small} \| (P_{\leq M+(N')^{1/4}} -
P_{\leq M})v \|_{Y^{s+1}} + \| (P_{\leq M+(N')^{1/4}} - P_{\leq M}) \tilde v \|_{Y^{s+1}} \leq C
(N')^{-\sigma}.
\end{equation}
Fix this $M$.  We split
$$ v = v_{lo} + v_{med} + v_{hi}$$
where
$$ v_{lo} := P_{\leq M} v; \quad v_{med} = (P_{\leq M+(N')^{1/4}} - P_{\leq M})v; \quad v_{hi} :=
(1 - P_{\leq M+(N')^{1/4}}) v.$$
Thus from \eqref{vtv-est}, \eqref{med-small} we have \be{vvv-est} \| v_{lo} \|_{Y^{s+1}}, \|
v_{hi} \|_{Y^{s+1}} \leq C; \quad \| v_{med} \|_{Y^{s+1}} \leq C(N')^{-\sigma}.
\end{equation}

Applying $P_{\leq M}$ to \eqref{mkdv-drift} and using Lemma \ref{form-lemma}, we see that $v_{lo}$
obeys the equation
$$
( \partial_t + \partial_{xxx}) v_{lo} = P_{\leq M} F_0(v, v, v) +  P_{\leq M} F_{\neq 0}(v, v,
v).$$ From the definition \eqref{f0-def} of the resonant operator $F_0$, we see that
$$ P_{\leq M} F_0(v,v,v) = F_0(v_{lo}, v_{lo}, v_{lo}).$$
The situation for $F_{\neq 0}$ is more complicated as this nonlinearity will mix $v_{lo}$,
$v_{med}$, $v_{hi}$ together.  Define an \emph{error term} to be any quantity with a $Z^{s+1}$
norm of $O((N')^{-\sigma})$.  From \eqref{vvv-est} and \eqref{fn0} we see that any term in
$F_{\neq 0}(v,v,v)$ involving $v_{med}$ is an error term.

Now let us consider the terms which involve $v_{hi}$.  A typical term is
$$ P_{\leq M} F_{\neq 0}(v_{lo}, v_{lo}, v_{hi}).$$
We can dyadically decompose this as
$$ \sum_{N_0, N_1, N_2, N_3} P_{N_0} P_{\leq M} F_{\neq 0}( P_{N_1} v_{lo}, P_{N_2} v_{lo},
 P_{N_3} v_{hi} ).$$
Such a term can be estimated using the frequency separation between $v_{lo}$ and $v_{hi}$: for the summand
to be nonzero, we need $N_1, N_2 \leq M$, and $N_3 \geq M + (N')^{\frac{1}{2}}$.  Using the notation in
the definition of $F_{\neq 0}$ \eqref{fn0-def}, we also need $|k_1 + k_2 + k_3| \sim N_0 \leq M$, hence
we must clearly also have
$N_{tenor} \gtrsim (N')^{1/4}$. 
From our non-resonant estimate
\eqref{fn0-improv}, the bounds
\eqref{vvv-est} above, and a summation of the dyadic indices $N_j$ (conceding some powers of $\log N'$ if
necessary) we thus see that this term is an error term.  A similar argument shows that any other
term involving $v_{hi}$ will also be an error term.  Thus we see that $v_{lo}$ obeys the equation
\begin{equation} \label{wpstuff}
( \partial_t + \partial_{xxx}) v_{lo} = F_0(v_{lo}, v_{lo}, v_{lo}) +  P_{\leq M} F_{\neq
0}(v_{lo}, v_{lo}, v_{lo}) + \hbox{error terms}.
\end{equation} 
By similar reasoning, the function $\tilde
v_{lo} := P_{\leq M} \tilde v$ also obeys the same equation (but with slightly different error
terms, of course).  Since $\tilde v_{lo}(0) = v_{lo}(0)$, we thus see from the standard local
well-posedness theory\footnote{A rough sketch of what we have in mind here is:  write $G$ for that portion 
of the nonlinearity on the right side of \eqref{wpstuff} not involving the error terms, and note
$$  \tilde v_{lo} - v_{lo}  = \int_0^t e^{i(t - \tau)\xi^3} \left(G(\tilde v_{lo}) - G(v_{lo}) + \text{error terms}
\right) d \tau.$$ Writing $G(\tilde v_{lo}) - G( v_{lo}) = \int_0^1 DG(\theta \tilde v_{lo} + (1 - \theta)
v_{lo}) (\tilde v_{lo} -  v_{lo})) d \theta$, we use \eqref{energy-est}, \eqref{f0}, \eqref{fn0}, 
and the fact that by scaling, we may assume that the data for $v_{lo}, \tilde v_{lo}$ are small 
in $Y^{s+1}$ to conclude \eqref{boundboundbound}.}  that
\begin{equation} \label{boundboundbound} \| \tilde v_{lo} - v_{lo} \|_{Y^{s+1}} \leq C (N')^{-\sigma}
\end{equation}
which by \eqref{sup-est} implies \eqref{v-v} as desired.  This proves Theorem
\ref{low-perturb-kdv}.

\section{Proof of Theorem \ref{kdv-approx}:  $BKdV$ approximates $KdV$ at low frequencies.}
\label{kdv-approx-sec}

We now prove the more difficult of our KdV approximation theorems, namely Theorem
\ref{kdv-approx}.  The proof here is definitely  in the same spirit as that of Theorem
\ref{low-perturb-kdv}, in that we show two flows remain close by showing that their mKdV analogues
remain close.  However, the proof will be more complicated since  one of the flows being
studied is $S_{B KdV}$ (see \eqref{b-pn-kdv}), and the standard Miura transform $\M$
defined by \eqref{miura-def} seems an inappropriate tool with which to pull the $S_{B KdV}$ flow
back to an $mKdV$-type evolution, as it introduces a $v_x^2$ type nonlinearity on the right side
of \eqref{mkdv-drift} which is too rough for us to estimate.  Instead, we introduce
a {\em{modified}} Miura transform $\M_B$.  This strategy is illustrated in \eqref{figure},
where we have written $S_{BmKdV}$  
for the flow which intertwines  $\M_B$ and $B KdV$ in the sense that
\begin{align*}
\M_B \circ S_{BmKdV}(t) \circ \M_B^{-1} & \equiv S_{B KdV}(t).
\end{align*}
\begin{equation}
\label{figure}
\begin{CD}
 v_0 @>{S_{mKdV}(t)}>>  v(t) \\
@V{\M}VV @VV{\M}V \\
u_0 @>{S_{KdV}(t)}>>  u(t) \\
u_0 @>>{S_{B KdV}(t)}>  \tilde u(t) \\
@A{\M_B}AA @AA{\M_B}A \\
\tilde v_0 @>{S_{BmKdV}(t)}>> \tilde v(t)
\end{CD}
\end{equation}
We can summarize the proof of Theorem \ref{kdv-approx} (using the same notation as in 
\eqref{figure}, which will be
defined momentarily!) by saying that $u(t), \tilde u(t)$ are shown to be close at
low frequencies by showing that
$\tilde{v}(t), v(t)$ are likewise close.

We now turn to the details.  Fix $s \geq -1/2$, $T > 0$, $N \gg 1$, $B$, and $u_0 \in H^s_0$; our
implicit constants may depend on $s$, $T$, and $\| u_0 \|_{H^s_0}$. We work exclusively in the time interval $[-T,T]$.

Let $\tilde u(t) := S_{B KdV}(t) u_0$
denote the evolution of the flow \eqref{b-pn-kdv}.  Our task is to show that \be{kdv-approx-targ}
 \sup_{|t| \leq T} \| P_{\leq N^{1/2}}(S_{KdV}(t) u_0 - \tilde u(t)) \|_{H^s_0} \lesssim N^{-\sigma}.
\end{equation}

We first claim (in analogy with \eqref{kdv-bound}) the bound \be{b-kdv-bound} \sup_{|t| \leq T} \|
\tilde u(t) \|_{H^s_0} \lesssim 1,
\end{equation}
if $N$ is large enough. This bound is achieved by a repetition of the arguments in \cite{ckstt:2}.
As it is somewhat technical and uses techniques different from those elsewhere in this paper
(notably the ``$I$-method''), we defer the proof of \eqref{b-kdv-bound} to an Appendix.

We may assume from \eqref{b-kdv-bound} and the local well-posedness theory\footnote{The well-posedness 
theory for KdV from \cite{kpv:kdv} can be applied without substantial change to the BKdV 
equation \eqref{b-pn-kdv}.  The presence of the multiplier $B$ on the right hand side presents
no difficulty.} that $u_0$, and hence
$\tilde u$,  is smooth.   

The Miura transform \eqref{miura-def} intertwines the KdV flow with the (renormalized) mKdV flow \eqref{mkdv-drift},
\eqref{F-def}.  We seek a similar
transform to intertwine the KdV-like flow $S_{B KdV}$ with an mKdV-like flow.  It turns out
that the correct transform to use is given by
\begin{align} \label{star28}
 \M_B \tilde v & := \v_x + B(1-P_0) (\v^2) = \v_x + B(\v^2) - P_0(\v^2),
 \end{align}
where of course the multiplier $B$ here is that which appears in the flow \eqref{b-pn-kdv} above.

As with $\M$, the operator $\M_B$ is a locally Lipschitz map from $H^{s+1}_0$ to $H^s_0$.  We now
address the question of invertibility of $\M_B$.

Let $\v$ be a function bounded in $H^{s+1}_0$.  We first look at the derivative operator
$\M'_B(\v)$, defined by
$$ \M'_B(\v) f := f_x + 2 B(1-P_0)(\v f).$$

\begin{lemma}\label{d-invert}  Fix $\v \in H^{s+1}_0, s \geq - \frac{1}{2}$, and allow the
 implicit constants in this Lemma to depend on
$\|\v\|_{H^{s+1}_0}$.
If $N$ is sufficiently large, then $\M'_B(\v)$ is invertible from $H^s_0$ to $H^{s+1}_0$, in the
sense that
$$ \| \M'_B(\v)^{-1} f \|_{H^{s+1}_0} \lesssim \| f \|_{H^s_0}$$
for all (smooth) $f$.
\end{lemma}

\begin{proof}
Recall from the proof of Theorem \ref{miura-invert} that we have the bound \be{m-inv}
 \| \M'(\v)^{-1} f \|_{H^{s+1}_0} \lesssim \| f \|_{H^s_0}.
\end{equation}
We proved this for $s= - 1/2$ but it is easy to see the same argument works for
$s > - 1/2$.  From the resolvent identity 
$$ O^{-1} = A^{-1} (1 - (A-O)A^{-1} )^{-1}$$
it thus suffices to show that the operator
$$(\M'_B(\v) - \M'(\v)) \M'(\v)^{-1}$$
is a contraction on $H^s_0$.  Applying \eqref{m-inv} again, it thus suffices to show the bound
$$ \| \M'_B(\v) f - \M'(\v) f \|_{H^s_0} \ll \| f \|_{H^{s+1}_0}.$$
But the left-hand side is just
$$ \| 2 (1-B)(\v f) \|_{H^s_0} \lesssim N^{-\sigma} \| \v f \|_{H^{s+\sigma}} \lesssim N^{-\sigma}
\| \v \|_{H^{s+1}_0}
\| f \|_{H^{s+1}_0} \lesssim N^{-\sigma} \|f\|_{H^{s+1}_0}$$ by \eqref{sob-mult} for some $\sigma
> 0$, and the claim follows if $N$ is sufficiently large.
\end{proof}

\begin{corollary}\label{mb-invert}  Let $R > 0, s \geq -\frac{1}{2}$.  If $N$ is large enough
depending on $R$, then there
is a map
$\M_B^{-1}$ defined on the ball $\B^\infty(0; R) := \{ \tilde u \in H^s_0: \| \u \|_{H^s_0} \leq R \}$
which inverts $\M_B$ and is a Lipschitz map from $\B^\infty(0;R)$ to $H^{s+1}_0$.
\end{corollary}

Remark:  Recall that $M_B$ depends on $N$ through the definition of 
$B$ (see \eqref{b-pn-kdv}).

\begin{proof}
Fix $R$; implicit constants are allowed to depend on $R$.

Let $\u \in \B^\infty(0; R)$.  To define $\M_B^{-1}$ at $\u$ we of course have to solve the equation
$$ \M_B \v = \u.$$
From Theorem \ref{miura-invert} we can find a $\vapr$, bounded in $H^{s+1}_0$, such that
$$ \M \vapr = \u.$$
We now apply the ansatz $\v = \vapr + \w$.  One easily checks, using \eqref{star28}, that $\w$ verifies the difference
equation
$$ \w_x + B(1-P_0)(2\vapr \w + \w^2) = (1-B)(\vapr^2)$$
or equivalently
$$ \w = \M'_B(\vapr)^{-1} (1-B)(\vapr^2) - \M'_B(\vapr)^{-1} B(1-P_0)(\w^2).$$
Since $\vapr$ is bounded in $H^{s+1}$ we see from Lemma \ref{d-invert} and \eqref{sob-mult} that
$$ \| \M'_B(\vapr)^{-1} (1-B)(\vapr^2) \|_{H^{s+1}} \lesssim N^{-\sigma}.$$
A contraction mapping argument again using Lemma \ref{d-invert} and \eqref{sob-mult} thus shows
that a solution $\w$ to the above difference equation exists and obeys the bound
$$ \| \w \|_{H^{s+1}} \lesssim N^{-\sigma}$$
if $N$ is sufficiently large.  In particular we see that $\M_B^{-1}$ exists at $\u$ and that
$\M_B^{-1}$ is bounded on $H^{s}_0$.

The Lipschitz bound now follows from Lemma \ref{d-invert} and the inverse function theorem, since
$\M_B$ is a smooth map from $H^{s+1}_0$ to $H^s_0$.  (Equivalently, one can use contraction
mapping arguments similar to the one above to show that $\M_B^{-1}$ is uniformly Lipschitz on very
small neighbourhoods of $\u$, and hence on the whole ball $\B^\infty(0; R)$).
\end{proof}

Thus if $N$ is large enough, the above corollary and \eqref{b-kdv-bound} let us write
\begin{align}
\label{v-def}
\v(t) \equiv \M_B^{-1} \u(t)
\end{align}
and conclude also that,
\be{v-bound-s1} \sup_{|t| \leq T} \| \v(t) \|_{H^{s+1}_0} \lesssim 1.
\end{equation}
From the Leibnitz rule we see that \bas
\u_t &= \M'_B(\v) \v_t \\
\u_x &= \M'_B(\v) \v_x = \v_{xx} + 2 B(\v \v_x)\\
\u_{xxx} &= \M'_B(\v) \v_{xxx} + 6 B (\v_x \v_{xx})\\
\u u_x &= (\v_x + B(\v^2) - P_0(\v^2)) \M'_B(\v) \v_x \\
&= \v_x \v_{xx} + 2 \v_x B(\v \v_x) + B(\v^2) \v_{xx} + 2 B(\v^2) B(\v \v_x) - \M'_B(\v) (P_0(\v^2) \v_x)
\end{align*}
where we have used the fact that $P_0(f f_x) = 0$ for any $f$.  Expanding out \eqref{b-pn-kdv} and
canceling the two terms of $6 B(\v_x \v_{xx})$ which appear, we  obtain
$$ \M'_B(\v)( \v_t + \v_{xxx}) = 6 B( 2 \v_x B(\v \v_x) + B(\v^2) \v_{xx} + 2 B(\v^2) B(\v \v_x) ) - B \M'_B(\v)
(6 P_0(\v^2) \v_x).$$
The first term of the right-hand side is roughly $\M'_B(\v)(6 B(B(\v^2) \v_x))$.  Indeed, a computation
shows   
$$
\M'_B(\v)(6 B(B(\v^2) \v_x)) = 6 B( 2 \v_x B(\v \v_x) + B(\v^2) \v_{xx} ) + 12 B (1-P_0)(\v B(B(\v^2)
\v_x)).$$ Thus we have
$$ \M'_B(\v)( \v_t + \v_{xxx} - 6 B(B(\v^2) \v_x) + 6 B(P_0(\v^2) \v_x)  ) = 12 E_1 + 6 E_2$$
where the error terms $E_1$, $E_2$ are the ``commutator expressions'' \bas
E_1 &:= B( B(\v^2) B(\v \v_x) - (1-P_0)(\v B(B(\v^2) \v_x)) ) \\
E_2 &:= P_0(\v^2) [\M'_B(\v), B] \v_x
\end{align*}

Thus $\v$ obeys the equation \be{v-eq-error} \v_t + \v_{xxx} = 6 B( (B-P_0)(\v^2) \v_x ) +
\M'_B(\v)^{-1}( 12 E_1 + 6 E_2); \quad \v(0) = \v_0.
\end{equation}
We have written $S_{BmKdV}(t)$ in figure \eqref{figure} to represent this flow.
Since $\v$ is smooth, it is \emph{a priori} in the space $Y^{s+1}$ when restricted to the interval
$[-T,T]$.  We now seek to control the non-linear terms in \eqref{v-eq-error}.

If it were not for the error terms $E_1$, $E_2$, one could obtain bounds of the form
\be{vys-bound} \| \v \|_{Y^{s+1}} \lesssim 1
\end{equation}
from \eqref{v-bound-s1} and the local well-posedness theory for mKdV in \cite{ckstt:gKdV} (which
can easily handle the presence of the order 0 operator $B$).  To deal with the $E_1$, $E_2$ terms
we use the following estimate.

\begin{lemma}\label{v-eq-est}  We have
\be{me-error} \|\M'_B(\v)(t)^{-1} E_j \|_{Z^{s+1}} \leq C N^{-\sigma}.
\end{equation}
for $j=1,2$, and $t \in [-T,T]$.
\end{lemma}

\begin{proof}
By \eqref{crude} and Lemma \ref{d-invert} (using \eqref{v-bound-s1}, of course) it suffices to
show that \be{ej-bound}
 \| E_j \|_{L^\infty_t H^s_x} \lesssim N^{-\sigma}.
\end{equation}
We first prove this for $E_1$.  Observe that $B(\v^2) B(\v \v_x) = \partial_x \frac{1}{4} (B(\v^2))^2$
has mean zero, and so we can factor out a $(1-P_0)$, and reduce to showing that
$$
\| \w B(\v \v_x) - \v B(\w \v_x) \|_{H^s_x} \lesssim N^{-\sigma}$$ where we have used the shorthand $\w
:= B(\v^2)$.

By \eqref{sob-mult} we see that $\w$ is bounded in $H^{s+\sigma}_x$ for some $\sigma > 0$.  From
the identity
$$ \w B(\v \v_x) - \v B(\w \v_x) = \w [B,\v] \v_x - \v [B,\w] \v_x$$
and another application of \eqref{sob-mult}, we see that it suffices to show the commutator
estimate \be{comm} \| [B, f] g \|_{H^s_x} \lesssim N^{-\sigma/2} \| f \|_{H^{s+\sigma}_x} \| g
\|_{H^s_x}.
\end{equation}
Without loss of generality we may assume that $f$ and $g$ have non-negative Fourier transform.
Observe that
$$ \widehat{[B,f] g}(k) = \sum_{k_1 + k_2 = k} (b(k) - b(k_2)) \widehat f(k_1) \widehat g(k_2).$$
The quantity $b(k) - b(k_2)$ is clearly $O(1)$.  If $|k_1| \ll N$ then one also obtains a bound of
$O(|k_1|/N)$ by the mean-value theorem.  Thus we have a universal bound of
$$ |b(k) - b(k_2)| \lesssim |k_1|^{\sigma/2} N^{-\sigma/2}.$$
The commutator estimate then reduces to
$$ \| (|\partial_x|^{\sigma/2} f) g \|_{H^s_x} \lesssim  \| f \|_{H^{s+\sigma}_x} \| g \|_{H^s_x},$$
but this follows from \eqref{sob-mult}.

Now we prove \eqref{ej-bound} for $E_2$.  From \eqref{v-bound-s1} we see that $P_0(\v^2)$ is
bounded in time, so it suffices to show that
$$
\| [\M'_B(\v), B] \v_x \|_{L^\infty_t H^s_x} \lesssim N^{-\sigma}.$$ Since $[\partial_x, B] = 0$, we
have
$$ [\M'_B(\v), B] \v_x = B(1-P_0)(\v B \v_x) - B^2(1-P_0)(\v \v_x) = B(1-P_0) [\v,B] \v_x,$$
and the claim follows from \eqref{comm}.
\end{proof}

From this lemma and perturbation theory in the $Y^{s+1}$ spaces (using the local well-posedness
theory in \cite{ckstt:gKdV}) we thus obtain \eqref{vys-bound}.

We now repeat the argument from Section \ref{low-kdv-sec}. Recall the notation from figure \eqref{figure}
that $v(t) \equiv S_{mKdV}(t) v_0$.   From \eqref{kdv-bound},
\eqref{intertwining} and Theorem \ref{miura-invert} we see that $v(t)$ is uniformly bounded in
$H^{s+1}$.  From the local well-posedness theory for mKdV we thus have
$$ \|v \|_{Y^{s+1}} \lesssim 1.$$
From this and \eqref{vys-bound}, we may find an interval $[M, M + N^{1/4}] \subseteq [N^{1/2},
2N^{1/2}]$ such that
$$ \| (P_{\leq M+N^{1/4}} - P_{\leq M}) \v \|_{Y^{s+1}} +
\| (P_{\leq M+N^{1/4}} - P_{\leq M}) v \|_{Y^{s+1}} \lesssim N^{-\sigma}.$$ Fix this
$M$.  Set
$$\v_{lo} := P_{\leq M} \v \quad  \text{and} \quad    v_{lo}(t) := P_{\leq M} v.$$
By
arguing as in the previous section we see that $v_{lo}$ obeys the equation \be{vm-eq}
(\partial_t + \partial_{xxx})  v_{lo} = F_0( v_{lo},  v_{lo},  v_{lo}) +
F_{\neq 0}( v_{lo}, v_{lo},  v_{lo}) + \hbox{ error terms}
\end{equation}
where the error terms have a $Z^{s+1}$ norm of $O(N^{-\sigma})$. We now claim that $\v_{lo}$ obeys the
same equation (but with a different set of error terms, of course).  Assuming this claim for the moment,
note that $v_{lo}$ and $ \v_{lo}$ have the same initial data, so we obtain,
\begin{align} \label{nowfinished}
 \sup_{|t| \leq T} \| v_{lo}(t) - \tilde v_{lo}(t) \|_{H^{s+1}_0} & \lesssim N^{-\sigma} \end{align}
by perturbation theory. 
The bound \eqref{nowfinished}  implies our goal
\eqref{kdv-approx-targ} relatively quickly: apply the Miura transform $\M$ (see \eqref{miura-def})
to the difference on the left side of \eqref{nowfinished}, and use the commutator bound
\eqref{commutator}, the fact that $P_{\leq M} \M \equiv P_{\leq M} \M_B$, and $M \geq
N^{\frac{1}{2}}$ to conclude that,
\begin{align}
N^{-\sigma} & \gtrsim \| P_{\leq M} \M_B \v(t) - P_{\leq M} \M  v(t) \|_{H^s_0} \nonumber
\\
& \gtrsim \|P_{\leq N^{\frac{1}{2}}} \u(t) - P_{\leq N^{\frac{1}{2}}} u \|_{H^s_0},
\label{almostthere}
\end{align}
as desired (see \eqref{kdv-approx-targ}). 

It remains to show that $\v_{lo}$ verifies \eqref{vm-eq}. Applying $P_{\leq M}$ to
\eqref{v-eq-error} and using Lemma \ref{v-eq-est} we have
$$
(\partial_t + \partial_{xxx}) \v_{lo} = 6 P_{\leq M}( (B-P_0)(\v^2) \v_x ) + error.$$ By repeating
the argument in Section \ref{low-kdv-sec} we have
$$ 6 P_{\leq M}( (1-P_0)(\v^2) \v_x ) = P_{\leq M}( F_0(\v,\v,\v) + F_{\neq 0}(\v,\v,\v) ) =
F_0(\v_{lo},\v_{lo},\v_{lo}) + F_{\neq 0}(\v_{lo},\v_{lo},\v_{lo}) + \hbox{error terms}.$$
Thus it will suffice to show that \be{v-leftover} P_{\leq M} ( (1-B)(\v^2) \v_x ) = \hbox{error
terms}.
\end{equation}
For a fixed time $t$, the spatial Fourier coefficient of the left-hand side at $(k,t)$ is
$$ \sum_{k=k_1+k_2+k_3} \chi_{[-M,M]}(k) (1-b(k_1+k_2)) \widehat \v(k_1,t) \widehat \v(k_2,t) ik_3
\widehat \v(k_3,t).$$
The summand vanishes unless $|k| \leq M \lesssim N^{1/2}$ and $|k_1+k_2| \gtrsim N$, which forces
$|k_3| \gtrsim N$.

First consider the contributions of the case when $(k_1+k_2)(k_2+k_3)(k_1+k_3) \neq 0$.  We now
apply \eqref{fn0-improv}.  By our previous discussion we have $N_0 \lesssim N^{1/2}$ and
$N_{soprano} \gtrsim N$, hence we see from \eqref{vys-bound} (writing things in terms of spacetime
Fourier transforms instead of spatial Fourier transforms, taking absolute values and discarding
the $1-b(k_1+k_2)$ factor) that this contribution is $error$.

It remains to consider the case when $(k_1+k_2)(k_2+k_3)(k_1+k_3) =
0$. 
 By the previous
discussion $k_1+k_2$ cannot be zero, while $|k_3|$ is much larger than $|k|$.  Thus the only two
cases are when $(k_1,k_2,k_3)$ is equal to $(k,-k_3,k_3)$ or $(-k_3, k, k_3)$, so by symmetry the
total contribution to the Fourier coefficient is
$$ 2 \chi_{[-M,M]}(k) \sum_{|k_3| \gtrsim N} ik_3 (1-b(k_3-k)) \widehat \v(k, t) \widehat \v(-k_3,t) \widehat \v(k_3,t).$$
Combining the $k_3$ term with the $-k_3$ term, this becomes
$$ 2 \chi_{[-M,M]}(k) \sum_{k_3 \gtrsim N} ik_3 (b(-k_3-k)-b(k_3-k)) \widehat \v(k, t) \widehat \v(-k_3,t)
\widehat \v(k_3,t).$$
By the mean-value theorem and the fact that $b$ is even, we have
$$ (b(-k_3-k)-b(k_3-k)) = O(|k|/N) = O(N^{-\sigma}).$$
Meanwhile, we have
$$ \sum_{k_3 \gtrsim N} |k_3| |\widehat \v(-k_3,t)| |\widehat \v(k_3,t)| \lesssim \| \v \|_{H^{s+1}_0}^2
 \lesssim 1.$$
Thus the above Fourier coefficient is $O(N^{-\sigma} |\widehat \v(k,t)|)$.  By \eqref{v-bound-s1} we
thus see that this contribution to \eqref{v-leftover} has an $L^\infty_t H^{s+1}_0$ norm of
$O(N^{-\sigma})$.  By \eqref{crude} we thus see that this contribution is $error$ as desired.
This completes the proof of \eqref{vm-eq} and hence \eqref{kdv-approx-targ}. This concludes the
proof of Theorem \ref{kdv-approx}.

\section{Proof of Theorem \ref{main}:  Symplectic nonsqueezing of KdV}
\label{finale-sec}

Let $N \gg 1$, and let $b$ be a symbol adapted to $[-N,N]$ which equals one on $[-N/2,N/2]$, and
let $B$ be the associated Fourier multiplier.  We begin by considering the modified Hamiltonian
$H_N$ on $P_{\leq N} H^{-1/2}_0(\T)$, defined by
$$ H_N(u) := \int_\T -\frac{1}{2} u_x^2 - (Bu)^3\ dx.$$
We compute the Hamiltonian flow on $P_{\leq N} H^{-1/2}_0$ corresponding to $H_N$.  Fix $u, v \in
H^{-1/2}_0$.  We see that \bas \frac{d}{d\eps} H_N(u + \eps v)|_{\eps=0}
=& \int_\T - u_x v_x - 3(Bu)^2 Bv\ dx  \\
=& \{ -u_{xxx} + 6 B( (Bu) (Bu_x)), v \}.
\end{align*}

Since $-u_{xxx} + 6 B( (Bu) (Bu_x) )$ is in $P_{\leq N} H^{-1/2}_0$, we conclude as in \eqref{secondd1},
\eqref{secondd2} that the Hamiltonian
flow of $H_N$ on $P_{\leq N} H^{-1/2}_0$ is given by
\be{blub} u_t + u_{xxx} = 6 B( (Bu) (Bu_x) );\quad
u(0) = u_0 \in P_{\leq N} H^{1/2}_0(\T).
\end{equation}
Let $S^{(N)}_{KdV}(t)$ denote the flow map associated to this equation; for each $t$, we observe
that $S^{(N)}_{KdV}(t)$ is thus a symplectomorphism on the finite-dimensional symplectic vector
space $P_{\leq N} H^{-1/2}_0$.  In particular, it obeys Theorem \ref{finite} (that is, we pick $S^{(N)}_{\text{Good!}} \equiv S^{(N)}_{KdV})$).   To conclude the
proof of Theorem \ref{main} it thus suffices
to show that the flow $S^{(N)}_{KdV}(t)$ obeys the
weak approximation property in Condition \ref{kdv-lowfreq}:

\begin{proposition}\label{kdv-lowfreq-correct}  Let $k_0 \in \Z^*$, $T > 0 $, $A > 0$, $0 < \eps \ll 1$.
 Then there exists an $N_0 = N_0(k_0,T,\eps,A) \gg |k_0|$ such that
$$ |k_0|^{-1/2} |\widehat{S_{KdV}(T) u_0}(k_0) - \widehat{S^{(N)}_{KdV}(T) u_0}(k_0)| \leq \eps$$
for all $N \geq N_0$ and all $u_0 \in \B^N(0,A)$ (see \eqref{BN-def} for the definition of this ball).
\end{proposition}

\begin{proof}
We make the transformation $w := Bu$, where $u$ solves \eqref{blub}.  Applying $B$ to \eqref{blub} we obtain
$$
w_t + w_{xxx} = 6 B^2(w w_x); \quad w(0) = B u_0$$ which is \eqref{b-pn-kdv} with $B$ replaced by
$B^2$.  Thus we have the intertwining relationship described by \eqref{hiccupcure} in the introduction 
to this paper,
\begin{equation*}
 B S^{(N)}_{KdV}(t) u_0 = S_{B^2 KdV}(t) B u_0.
 \end{equation*}
In particular, if $N_0 \gg |k_0|$, then $b(k_0) = 1$, so we have
\begin{equation}
\label{twostarry}
 \widehat{S^{(N)}_{KdV}(T) u_0}(k_0) = \widehat{(S_{B^2 KdV}(T) B u_0)}(k_0).
 \end{equation}

From Theorem
\ref{low-perturb-kdv} we have
\begin{equation}
\label{first}
 |k_0|^{-1/2} | \widehat{S_{KdV}(T) u_0}(k_0)  - \widehat{S_{KdV}(T) Bu_0}(k_0) | \lesssim N^{-\sigma}.
 \end{equation}
From Theorem \ref{kdv-approx} we have (if $N_0$ is large enough, $N_0 \gg k_0$)
\begin{equation}
\label{second}
 |k_0|^{-1/2} | \widehat{S_{KdV}(T) B u_0}(k_0) -   \widehat{S_{B^2 KdV}(T) B u_0}(k_0)|
\lesssim N^{-\sigma},
\end{equation}
where the implicit constants are allowed to depend on $T$ and $A$.  By \eqref{twostarry},
the second term on the left of \eqref{second} is the same as $\widehat{S^{(N)}_{KdV}(T)
u_0}(k_0)$.
Combining this observation with \eqref{first}, \eqref{second}, and the triangle inequality,
we obtain the desired claim, if $N_0$ is sufficiently large depending
on $k_0$, $T$, $\eps$, $A$.
\end{proof}

The proof of Theorem \ref{main} is now complete.

\section{Proof of Theorem \ref{kdv-lowfreq-uhoh}: $P_{\leq N}$$KdV$ does not approximate $KdV$}
\label{counter-sec}

Informally, the point of this section is that there is absolutely no slack in the bilinear
estimate \eqref{bil} at regularity $s=-1/2$ no matter what the frequencies of the various
functions are; see the examples in \cite{kpv:kdv}.  But to convert the examples for the bilinear
estimate to quantitative estimates of the KdV and truncated KdV flow - in particular, to establish
that the two flows differ as claimed in Theorem \ref{kdv-lowfreq-uhoh} - we must do some tedious
computation of iterates, which we detail below.

Fix $k_0$,$A$, $T$, for instance $T, A \sim 1$; our implicit constants in this section will be allowed to depend on these
parameters. 
Without loss of generality we may assume that $k_0 >
0$.  We let $0 < \sigma \ll 1$ be a small parameter depending on $k_0$, $A$, $T$ to be chosen
later.

Let $N \gg \sigma^{-100}$ be a large integer.  We consider the initial data
$$ u_0(x) := \sigma^3 \cos(k_0 x) + \sigma N^{1/2} \cos(N x).$$
Note that $u_0$ lies in $P_{\leq N} H^{-1/2}_0(\T)$ with norm $O(\sigma)$, and in particular we have $u_0
\in \B^N(0; A)$ if $\sigma \ll 1$ is sufficiently small.

Let $u$ and $u^{(N)}$ be the solutions to the KdV flow \eqref{kdv} and truncated KdV flow
\eqref{pn-kdv} respectively, with initial data $u(0) = u^{(N)}(0) = u_0$.  We shall show that, if
$\sigma$ is sufficiently small, \be{uun} |\widehat{u(T)}(k_0) - \widehat{u^{(N)}(T)}(k_0)| \sim
\sigma^5,
\end{equation}
which gives \eqref{fluctuate}.

To prove \eqref{uun} we need good approximations of $u$ and $u^{(N)}$.  To approximate $u$, we
look at the iterates $u^{[j]}$ for $j= 0,1,2,\ldots$ defined inductively by $u^{[-1]}(t,x)
\equiv 0$ and
\begin{align} \label{blahblah}
(\partial_t + \partial_{xxx}) u^{[j]} = \partial_x (3 (u^{[j-1]})^2); \quad u^{[j]}(0) = u_0.
\end{align}
From the contraction mapping arguments in \cite{kpv:kdv} (see also \cite{ckstt:gKdV}) we know that
the $u^{[j]}$ converge to $u$ in the $Y^{-1/2}_{[0,T]}$ norm; indeed each iterate is closer to $u$
by a factor of at least $O(\sigma)$ compared to the previous one\footnote{Strictly speaking, this
contraction mapping property was only proven for $T$ sufficiently small, but by subdividing
$[0,T]$ into a finite number of small intervals one can obtain the same contraction mapping for
arbitrary $T$ if $\sigma$ is sufficiently small depending on $T$.  This naive argument requires
$\sigma \ll e^{-CT}$ for some $C$; the more sophisticated scaling argument in \cite{ckstt:gKdV}
can improve this to $\sigma \ll T^{-1/3-}$, but we will not need this quantitative improvement for
our arguments here.}. A routine calculation yields
$$
u^{[0]}(t,x) = \sigma^3 \cos(k_0 x + k_0^3 t) + \sigma N^{1/2} \cos(Nx + N^3 t),$$ and thus
$$ \partial_x (3 (u^{[0]})^2) =
- \frac{3}{2} \sigma^4 N^{3/2} \sin((N+k_0)x + (N^3+k_0^3)t) - \frac{3}{2} \sigma^4 N^{3/2}
\sin((N-k_0)x + (N^3-k_0^3)t) + O_{Z}(\sigma^6)$$ where $O_Z(K)$ denotes a quantity with a
$Z^{-1/2}_{[0,T]}$ norm of $O(K)$ (note that we have used the hypothesis $N \gg \sigma^{-100}$ to
absorb several terms into this $O_Z(\sigma^6)$ error\footnote{For example, the term $\sigma^4
N^{\frac{1}{2}} k_0 \sin((N+k_0)x + (N^3+k_0^3)t)$ which appears when one calculates $\partial_x
(3 (u^{[0]})^2)$ is  $O_{Z}(\sigma^6)$, as the space-time Fourier transform of this term is
supported a distance approximately $N^2$ from the cubic $\tau = \xi^3$.  Hence when computing the
$Z^{-\frac{1}{2}}$ norm of this term, we get a factor of $N^{-1} \ll \sigma^{100}$ from the
denominator in the definition of this norm.}).

Observe that
$$ (\partial_t + \partial_{xxx}) \left(-\frac{1}{2} \sigma^4 N^{-1/2}
[\cos((N+k_0)x + (N^3+k_0^3)t) - \cos((N+k_0)x+ (N + k_0)^3 t)])\right)$$
$$
= - \frac{3}{2} \sigma^4 N^{3/2} k_0 \sin((N+k_0)x + (N^3+k_0^3)t) + O_Z(\sigma^6)$$ and
$$ (\partial_t + \partial_{xxx}) \left(\frac{1}{2} \sigma^4 N^{-1/2}
(\cos((N-k_0)x + (N^3-k_0^3)t) - \cos((N-k_0)x + (N - k_0)^3t))\right)$$
$$
= - \frac{3}{2} \sigma^4 N^{3/2} k_0 \sin((N-k_0)x + (N^3-k_0^3)t) + O_Z(\sigma^6).$$ Combining
this with the calculation of $\partial_x (3 u^{[0]})^2$ above and using \eqref{energy-est} we obtain
\bas u^{[1]}(t,x) = &u^{[0]}(t,x) \\
&+ \left(-\frac{1}{2} \sigma^4 N^{-1/2} k_0^{-1}
(\cos((N+k_0)x + (N^3+k_0^3)t) - \cos((N+k_0)x + (N + k_0)^3t))\right)\\
&+ \left(\frac{1}{2} \sigma^4 N^{-1/2} k_0^{-1}
(\cos((N-k_0)x + (N^3-k_0^3)t) - \cos((N-k_0)x + (N - k_0)^3 t))\right)\\
&+ O_Y(\sigma^6),
\end{align*}
where $O_Y(\sigma^6)$ denotes a quantity with a $Y^{-1/2}_{[0,T]}$ norm of $O(\sigma^6)$. In fact,
since the $\cos((N\pm k_0)x + (N \pm k_0)^3t)$ terms are already $O_Y(\sigma^6)$ we have
$$
u^{[1]}(t,x) = u^{[0]}(t,x) + \frac{1}{2} \sigma^4 N^{-1/2} k_0^{-1} \left(\cos((N-k_0)x +
(N^3-k_0^3)t) - \cos((N+k_0)x + (N^3+k_0^3)t)\right) + O_Y(\sigma^6).$$ Using \eqref{bil} to
handle any interaction with a factor of $\sigma^6$ or better, we obtain \be{0-1-approx}
\partial_x (3 (u^{[1]})^2) =
\partial_x (3 (u^{[0]})^2) + O_Z(\sigma^6).
\end{equation}
Note that there are two additional, potentially disruptive terms of the form $\pm \frac{3}{2}
\sigma^5 \sin(k_0 x + k_0^3 t)$ which appear in the expansion of $\partial_x (3 (u^{[1]})^2)$, but
they have opposite signs and so cancel\footnote{This special cancellation seems to be what
distinguishes the KdV flow \eqref{kdv} from superficially similar flows such as \eqref{pn-kdv},
and is crucial to obtaining our high-frequency and low-frequency approximation results for this
flow.  It is instructive to see this cancellation via the renormalized mKdV flow
\eqref{mkdv-drift} by computing iterates for mKdV and then applying the Miura transform to those
iterates.} each other. From \eqref{0-1-approx} and \eqref{energy-est} we have
$$ u^{[2]} = u^{[1]} + O_Y(\sigma^6).$$
From the contraction mapping property of the iteration map we thus have
$$ u = u^{[1]} + O_Y(\sigma^6).$$
In particular we see that \be{u1-coefficient} \widehat{u(T)}(k_0) = \widehat{u^{[1]}(T)}(k_0) +
O(\sigma^6) = \widehat{u^{[0]}(T)}(k_0) + O(\sigma^6).
\end{equation}

Now we approximate $u^{(N)}$.  To do this we construct iterates $\tilde u^{[j]}$, $j=0,1,2,\ldots$
for the truncated equation by setting $\tilde u^{[0]} := u^{[0]}$ and
$$ (\partial_t + \partial_{xxx}) \tilde u^{[j]} = P_{\leq N} \partial_x (3 (\tilde u^{[j-1]})^2);
 \quad \tilde u^{[j]}(0) = u_0.$$
By a variant of the local well-posedness theory from \cite{kpv:kdv} (and \cite{ckstt:gKdV}) we
know that $\tilde u^{[j]}$ will converge to $u^{(N)}$ in the $Y$ norm. By reviewing the
computation of $u^{[1]}(t,x)$, but now bearing in mind the presence of the projection $P_N$, we
obtain for the first iterate,
\begin{align*}
\tilde u^{[1]}(t,x) & = u^{[0]}(t,x) + \frac{1}{2} \sigma^4 N^{-\frac{1}{2}}
\cos((N-k_0)x + (N^3-k_0^3)t)) + O_Y(\sigma^6) \\
& = u^{[1]}(x,t) + O_{Y}(\sigma^6).
\end{align*}
 Comparing this with  the formula for $u^{[1]}$
above, we note that the Fourier modes at $\pm(N+k_0)$ are not present here. As a consequence, the
analog of \eqref{0-1-approx} reads,
$$
\partial_x (3 (\tilde u^{[1]})^2) =
\partial_x (3 (u^{[0]})^2) + \frac{3}{2} \sigma^5 \sin(k_0 x+ k_0^3 t) + O_Z(\sigma^6).
$$
Since $(\partial_t + \partial^3_x) (t \sin(k_0 x+ k_0^3 t)) = \sin(k_0 x+ k_0^3 t)$, we can write,
$$ \tilde u^{[2]} = \tilde u^{[1]} + \frac{3}{2} \sigma^5 t \sin(k_0 x+k_0^3 t) + O_Y(\sigma^6).$$
We can easily check then that,
\begin{align*}
\partial_x (3 (\tilde u^{[2]})^2)&  =
\partial_x (3 (u^{[1]})^2) + O_Z(\sigma^6),
\end{align*}
hence $\tilde u^{[3]} \, = \, \tilde{u}^{[2]} + O_Y(\sigma^6)$,
 which by the contraction mapping property implies
that
$$ u^{(N)} = \tilde u^{[2]} + O_Y(\sigma^6).$$
In particular we see that
$$
\widehat{u^{(N)}(T)}(k_0) = \widehat{\tilde u^{[2]}(T)}(k_0) + O(\sigma^6) =
\widehat{u^{[1]}(T)}(k_0) - \frac{3}{2} i T \sigma^5  e^{i k_0^3 T} + O(\sigma^6).
$$
Comparing this with \eqref{u1-coefficient} we obtain \eqref{uun} as desired.  This proves Theorem
\ref{kdv-lowfreq-uhoh}.

\section{ Appendix.  Proof of \eqref{b-kdv-bound}:  $H^s$ bound for the $B$KdV flow}
\label{appendix-sec}

We now prove the bound \eqref{b-kdv-bound} for $H^s_0$ solutions to the KdV-like equation
$$ u_t + u_{xxx} = 6 B(u u_x); \quad u(0) = u_0$$
with $\| u_0 \|_{H^s_0} \lesssim 1$; this bound is needed to complete the proof of Theorem
\ref{kdv-approx} and hence Theorem \ref{main}.

If $s \geq 0$ then this bound follows from $L^2$ conservation and standard persistence of
regularity theory (see e.g. \cite{borg:xsb}), so we shall assume that $-1/2 \leq s < 0$.

To do so, let us first review (from \cite{ckstt:2}) how the corresponding bound \eqref{kdv-bound}
was proven for the KdV flow
$$ u_t + u_{xxx} = 6 u u_x; \quad u(0) = u_0.$$

\subsection{Review of proof of $H^s$ bound for KdV \eqref{kdv-bound}}\label{kdv-bound-sec}

The idea is to modify the conserved $L^2$ norm $\int u^2$ to something resembling the $H^s$ norm
and which is still approximately conserved.  To do this it is convenient to introduce some
notation for multilinear forms.

If $n \geq 2$ is an integer, we define a \emph{(spatial) $n$-multiplier} to be any function
$M_n(k_1, \ldots, k_n)$ on the (discrete) hyperplane
$$ \Gamma_n := \{ (k_1, \ldots, k_n) \in \Z_*^n: k_1 + \ldots +
k_n = 0 \}.$$

If $M_n$ is a $n$-multiplier and $u_1, \ldots, u_n$ are functions on $\R/2\pi\Z$, we define the
\emph{$n$-linear functional} $\Lambda_n(M_n;u_1, \ldots, u_n)$ by
$$ \Lambda_n(M_n;f_1, \ldots, f_n)
:= \sum_{(k_1, \ldots, k_n) \in \Gamma_n} M_n(k_1, \ldots, k_n) \prod_{j=1}^n \widehat f_j(k_j).$$ We
adopt the notation
$$ \Lambda_n(M_n;u) := \Lambda_n(M_n; u, \ldots, u).$$
Observe that $\Lambda_n(M_n;f)$ is invariant under permutations of the $k_j$ indices.  In
particular we have
$$ \Lambda_n(M_n;u) = \Lambda_n([M_n]_{sym}; u)$$
where \be{sym} [M_n]_{sym}(k) := \frac{1}{n!} \sum_{\sigma \in S_n} M_n(\sigma(k))
\end{equation}
is the symmetrization of $M_n$.

Thus, for instance, we have $\int u^2 = 2\pi \Lambda_2(1; u)$, and more generally $\| u
\|_{H^s_0}^2 = 2 \pi \Lambda_2(|k_1|^s |k_2|^s; u) = 2 \pi \Lambda_2(|k_1|^{2s}; u)$ for $u \in
H^s_0$.

Now suppose that $u$ obeys the KdV evolution \eqref{kdv}, and $M_n$ is a symmetric multiplier.
Then we have the differentiation law \be{diff} \frac{d}{dt} \Lambda_n(M_n; u(t)) = \Lambda_n(M_n
\alpha_n; u(t)) - 3in \Lambda_{n+1}(M_n(k_1, \ldots, k_{n-1}, k_n+k_{n+1}) (k_n + k_{n+1}); u(t))
\end{equation}
where
$$ \alpha_n := k_1^3 + \ldots + k_n^3.$$
(see \cite{ckstt:2}).  Thus for instance we have \bas
 \frac{d}{dt} \Lambda_2(1; u(t)) &= \Lambda_2(\alpha_2; u(t))
- 6i \Lambda_3(k_2 + k_3; u(t)) \\
&= \Lambda_2(k_1^3 + k_2^3; u(t)) - 4i \Lambda_3(k_1+k_2+k_3; u(t)) \\
&= 0 - 0,
\end{align*}
demonstrating the conservation of the $L^2$ norm.

Henceforth we shall omit the $u(t)$ from the $\Lambda_n$ notation for brevity.  We also adopt the
convenient notation that $k_{ij} := k_i + k_j$, etc., thus for instance $k_{145} = k_1 + k_4 +
k_5$.  Also we write $m_i := m(k_i)$, $m_{ij} := m(k_{ij})$, etc, and $N_i$ for $|k_i|$, $N_{ij}$
for $|k_{ij}|$, etc.

Let $A \gg 1$ be a large number to be chosen later\footnote{ Note
that the quantity $A$ here represents what was called $N$ in \cite{ckstt:2}, a notational change necessary
since in the present paper $N$ represents something else.}, and let $m(k)$ be a multiplier which equals 1
on $[-A,A]$, equals $(|k|/A)^s$ for $|k| \geq 2A$, and is real, even, and smooth in between.  We
denote the corresponding Fourier multiplier by $I$:
$$ \widehat{Iu}(k) := m(k) \widehat u(k),$$
thus $I$ acts like the identity on frequencies $\leq A$ and is smoothing on frequencies $\gtrsim
A$.  We define the modified energy $E_2(t)$ by
$$ E_2(t) := \Lambda_2(m_1 m_2),$$
then one can verify that
$$ \| u(t)  \|_{H^s_0}^2 \lesssim E_2(t) \lesssim A^{-2s} \| u(t) \|_{H^s_0}^2.$$
From \eqref{diff}, \eqref{sym} and the fact that $\alpha_2 = 0$ we have \bas
\frac{d}{dt} E_2(t) &= -6i \Lambda_3(m_1 m_{23} k_{23}) \\
&= 6i \Lambda_3(m_1^2 k_1) \\
&= \Lambda_3(M_3)
\end{align*}
where $M_3$ is the 3-multiplier
$$ M_3 := 2i(m_1^2 k_1 + m_2^2 k_2 + m_3^2 k_3).$$
Now define the modified energy $E_3(t)$ by
$$ E_3(t) := E_2(t) + \Lambda_3(\sigma_3)$$
where $\sigma_3(k_1,k_2,k_3)$ is the 3-multiplier
$$ \sigma_3 := - M_3/\alpha_3.$$
This multiplier may appear to be singular at first glance, but we observe that \be{alpha-3}
\alpha_3 = k_1^3+k_2^3+k_3^3 = 3 k_1 k_2 k_3
\end{equation}
and that $M_3$ vanishes whenever $k_1 k_2 k_3=0$. Then by \eqref{diff}, \eqref{sym} we have \bas
\frac{d}{dt} E_3(t) &= \Lambda_3(M_3) + \Lambda_3(\sigma_3 \alpha_3) -
-9i \Lambda_4(\sigma_3(k_1,k_2,k_{34}) k_{34})\\
&= \Lambda_4(M_4)
\end{align*}
where $M_4$ is the 4-multiplier
$$ M_4 := -9i [\sigma_3(k_1,k_2,k_{34}) k_{34}]_{sym}.$$
Now define the modified energy $E_4(t)$ by
$$ E_4(t) := E_3(t) + \Lambda_4(\sigma_4)$$
where $\sigma_4k_1,k_2,k_3,k_4)$ is the 4-multiplier
$$ \sigma_4:= - M_4/\alpha_4$$
This multiplier may appear to be singular at first glance, but we observe that \be{alpha-4}
\alpha_4 = k_1^3+k_2^3+k_3^3+k_4^3 = 3 k_{12} k_{13} k_{14}
\end{equation}
(cf. \eqref{fund}), and one can check that $M_4$ vanishes when $k_{12} k_{13} k_{14} = 0$.  Then
as before we have that \be{e4-deriv} \frac{d}{dt} E_4(t) = \Lambda_5(M_5)
\end{equation}
where
$$ M_5 := -12i [\sigma_4(k_1,k_2,k_3,k_{45}) k_{45}]_{sym}.$$
We could continue this procedure indefinitely, but $E_4$ will turn out to be a suitable almost
conserved quantity for our purposes.  In \cite{ckstt:2} it was shown (by Gagliardo-Nirenberg type
arguments) that $E_4$ is bounded if and only if $\| u \|_{H^s_0}$ is bounded, so to obtain
\eqref{kdv-bound} it suffices to control $E_4(t)$.  In light of \eqref{e4-deriv} it will suffice
to control $M_5$.  The key lemma here was the following:

\begin{lemma}\label{m5-bound}\cite{ckstt:2} Let $k_1,k_2,k_3,k_4,k_5$ be real numbers (not necessarily
integer) such that $k_{12345} = 0$.  Then $M_5(k_1,\ldots,k_5)$ vanishes when $N_1, \ldots, N_5 \ll A$.
 In all other cases we have the bound
$$ |M_5(k_1, \ldots, k_5)| \lesssim [\frac{m^2(N_{*45}) N_{45}}{(A+N_1)(A+N_2)(A+N_3)(A+N_{45})}]_{sym}$$
where
$$ N_{*45} = \min(N_1, N_2, N_3, N_{45}, N_{12}, N_{13}, N_{14}).$$
\end{lemma}

With this bound and some multilinear $Y^s$ estimates\footnote{Strictly speaking, in order to
handle large data, these estimates had to take place in the large period setting
$\R/2\pi\lambda\Z$, as one would need to rescale large data to be small.  This causes some
unpleasant technical complications in the arguments, and in particular this is why the $k_j$ in
the above lemma need to be real (or lie in $\Z_*/\lambda$) rather than integer.  See
\cite{ckstt:2}, \cite{ckstt:gKdV} for more details.  In this paper we will ignore the large period
issue as it does not cause any essential change to the
argument.}, 
bounds on the growth of
$E_4(t)$ was obtained.  In particular, if $E_4(T)$ was small for some time $T$, it was possible to
obtain the bound $E_4(T+\delta) = E_4(T) + O(A^{-\frac{5}{2}-})$ for some small time $\delta \sim
1$.  Iterating this and using a rescaling argument one could obtain \eqref{kdv-bound} for all $s
\geq -1/2$ (after choosing $A$ appropriately depending on $\| u_0 \|_{H^s_0}$ and $T$).  See
\cite{ckstt:2} for details.

\subsection{Adapting the argument to the $B$KdV flow}

We now adapt the above argument to the flow \eqref{b-pn-kdv}.  The main difference will be the
appearance of various quantities of the form $b(k_i)$, $b(k_{ij})$, etc., however these factors
will play essentially no role in the argument.  Accordingly, we write $b_i$ for $b(k_i)$, etc.  We
shall assume that the frequency parameter $N$ corresponding to $b$ is much larger than the
frequency parameter $A$ corresponding to $m$.

Suppose $\tilde u$ solves \eqref{b-pn-kdv}.  Then \eqref{diff} now becomes
\be{diff-mod} \frac{d}{dt} \Lambda_n(M_n; \tilde u(t)) = \Lambda_n(M_n \alpha_n; \tilde u(t)) - 3in
\Lambda_{n+1}(M_n(k_1, \ldots, k_{n-1}, k_n+k_{n+1}) b(k_n + k_{n+1}) (k_n + k_{n+1}); \tilde u(t)).
\end{equation}
Again we define
$$ E_2(t) := \Lambda_2(m_1 m_2),$$
then one can verify that
$$
\frac{d}{dt} E_2(t) = \Lambda_3(M_3)$$ where $M_3$ is the 3-multiplier
$$ M_3 := 2i(f_1 + f_2 + f_3)$$
and $f(k) := m^2(k) b(k) k$.  Observe that $f$ is an odd function with $f'(k) = O(m(k))$ and
$f''(k) = O(m(k)/(A+|k|))$ for all $k$.

We observe the following bounds on $M_3$.

\begin{lemma}\label{M3-bound}  If $N_1, N_2, N_3 \ll A$, then $M_3 = 0$.  Otherwise, we have
$$ |M_3| \lesssim \max(m_1^2,m_2^2,m_3^2) \min(N_1,N_2,N_3).$$
\end{lemma}

\begin{proof}\cite{ckstt:2}
When $N_1, N_2, N_3 \ll A$ then $f_i=k_i$ for $i=1,2,3$, and the claim is clear.  Otherwise, we
use symmetry to assume that $N_1 \sim N_2 \gtrsim N_3$.  But then the mean-value theorem and the
above bounds on $f$ gives
$$ f_2 = - f_{13} = - f_1 + O(m^2_1 N_3),$$
and the claim easily follows.
\end{proof}

Now define the modified energy $E_3(t)$ by
$$ E_3(t) := E_2(t) + \Lambda_3(\sigma_3)$$
where $\sigma_3(k_1,k_2,k_3)$ is the 3-multiplier
$$ \sigma_3 := - M_3/\alpha_3.$$

From Lemma \ref{M3-bound}, \eqref{alpha-3} we see that $\sigma_3$ vanishes when $\max(N_1, N_2,
N_3) \ll N$, and we have the bounds \be{a3-bound} |\sigma_3| \lesssim
\frac{\max(m_1^2,m_2^2,m_3^2)}{(N + \max(N_1,N_2,N_3))^2}
\end{equation}
otherwise (note that the two largest values of $N_j$ have to be comparable).

By \eqref{diff}, \eqref{sym} we have
$$ \frac{d}{dt} E_3(t) = \Lambda_4(M_4)$$
where $M_4$ is the 4-multiplier
$$ M_4 := -9i [\sigma_3(k_1,k_2,k_{34}) b_{34} k_{34}]_{sym}.$$
Now define the modified energy $E_4(t)$ by
$$ E_4(t) := E_3(t) + \Lambda_4(\sigma_4)$$
where $\sigma_4(k_1,k_2,k_3,k_4)$ is the 4-multiplier
$$ \sigma_4:= - M_4/\alpha_4.$$
Then as before we have that
$$
\frac{d}{dt} E_4(t) = \Lambda_5(M_5)
$$
where
$$ M_5 := -12i [\sigma_4(k_1,k_2,k_3,k_{45}) b_{45} k_{45}]_{sym}.$$

Our aim is to show that this new $M_5$ still verifies the bounds in Lemma \ref{m5-bound}; the rest
of the arguments in \cite{ckstt:2} will then give the desired bound \eqref{b-kdv-bound} (the
presence of the $B$ multiplier having no impact on the local well-posedness  theory).

From the definition of $\sigma_4$ and $M_5$, it will suffice to prove the following $M_4$ bound.

\begin{lemma}\label{m4-bound}  If $\max(N_1,N_2,N_3,N_4) \ll A$ then $M_4$ vanishes.  Otherwise, we
have
$$ |M_4| \lesssim \frac{|\alpha_4| m^2(N_*)}{(A+N_1)(A+N_2)(A+N_3)(A+N_4)},$$
where $N_* := \min(N_1,N_2,N_3,N_4,N_{12}, N_{13}, N_{14}).$
\end{lemma}

\begin{proof}
When $\max(N_1,N_2,N_3,N_4) \ll A$ then $\sigma_3(k_1,k_2,k_{34})$ and all of its symmetrizations
vanish, hence $M_4$ vanishes.  Now we assume that $\max(N_1,N_2,N_3,N_4) \gtrsim A$.  By symmetry
we may assume that $N_1 \gtrsim N_2 \gtrsim N_3 \gtrsim N_4$, thus $N_1 \sim N_2 \gtrsim A$.  From
\eqref{alpha-4} we have $|\alpha_4| \sim N_{13} N_{14} N_{34}$.

We divide into several cases depending on the relative sizes of $N_2, N_3, N_4$.

{\bf Case 1: $N_2 \gg N_3 \gg N_4$.}  In this case $|\alpha_4| \sim N_1^2 N_3$, thus we reduce to
showing
$$ |M_4| \lesssim \frac{ m(N_*)^2 }{ A + N_4 }.$$
But from Lemma \ref{M3-bound} we have
$$ |\sigma_3(k_a, k_b, k_{cd}) b_{cd} k_{cd}| \lesssim \frac{\min(m_a, m_b, m_{cd})^2}{A +
\max(N_a, N_b, N_{cd})} \lesssim m(N_*)^2 / (A+N_4)$$
as desired.

{\bf Case 2: $N_2 \sim N_3 \gg N_4$.}  In this case $|\alpha_4| \sim N_1^3$, thus we reduce to
showing
$$ |M_4| \lesssim \frac{ m(N_*)^2 }{ A + N_4 }.$$
One then proceeds as in Case 1.

{\bf Case 3: $N_2 \gg N_3 \sim N_4$.}  In this case $|\alpha_4| \sim N_1^2 N_{34}$, thus we reduce
to showing
$$ |M_4| \lesssim \frac{ m(N_*)^2 N_{34}}{ (A + N_3)^2 }.$$
From Lemma \ref{M3-bound} we have
$$ |\sigma_3(k_1,k_2,k_{34}) b_{34} k_{34}| \lesssim m(N_*)^2 N_{34} / (A + \max(N_1, N_2, N_{34}))^2$$
which is acceptable.  Similarly
$$ |\sigma_3(k_3,k_4,k_{12}) b_{12} k_{12}| \lesssim m(N_*)^2 N_{12} / (A + \max(N_3, N_4, N_{12}))^2$$
is acceptable since $N_{12} = N_{34}$.  It thus suffices to show that
$$
| \sigma_3(k_1,k_3,k_{24}) b_{24} k_{24} + \sigma_3(k_1,k_4,k_{23}) b_{23} k_{23} +
\sigma_3(k_2,k_3,k_{14}) b_{14} k_{14} + \sigma_3(k_2,k_4,k_{13}) b_{13} k_{13} | \lesssim
\frac{m(N_*)^2 N_{34}}{(A + N_3)^2}.$$

We expand out $\sigma_3$ using \eqref{alpha-3}, and replace $k_1$ by $-k_{234}$ throughout, and
reduce to showing
$$
| -\frac{b_{24} (f_3 + f_{24} - f_{234})}{k_{234} k_3} - \frac{b_{23} (f_4 + f_{23} -
f_{234})}{k_{234} k_4} + \frac{b_{23} (f_2 + f_3 - f_{23})}{k_2 k_3} + \frac{b_{24} (f_2 + f_4 -
f_{24})}{k_2 k_4} | \lesssim \frac{m(N_*)^2 N_{34}}{(A + N_3)^2}.$$

From the mean-value theorem we have $b_{23} = b_2 + O(N_3/N_2) = b_2 + O(N_3/(A+N_3))$.  Similarly
$b_{24} = b_2 + O(N_3/(A+N_3))$.  Let us consider the contribution of the $O(N_3/(A+N_1))$ errors.
It will suffice to show that
$$
-\frac{f_3 + f_{24} - f_{234}}{k_{234} k_3}  + \frac{f_2 + f_4 - f_{24}}{k_2 k_4}$$ and
$$
- \frac{f_4 + f_{23} - f_{234}}{k_{234} k_4} + \frac{f_2 + f_3 - f_{23}}{k_2 k_3} $$ are both
$O(\frac{m(N_*)^2 N_{34}}{N_3 (A+N_3)})$. By the $k_3 \leftrightarrow k_4$ symmetry it suffices to
estimate the former expression.  From the mean-value theorem we have
$$ \frac{1}{k_{234}k_3} = \frac{1}{(k_2+k_{34}) (-k_4+k_{34})} = -\frac{1}{k_2 k_4} + O(\frac{N_{34}}
{N_2 N_4^2}).$$
By Lemma \ref{M3-bound}, the contribution of the error term $ O(\frac{N_{34}}{N_2 N_4^2})$ is
bounded by
$$ m(N_*)^2 N_3 O(\frac{N_{34}}{N_2 N_4^2}),$$
which is acceptable.  Thus it suffices to show that
$$
\frac{f_3 + f_{24} - f_{234}}{k_2 k_4}  + \frac{f_2 + f_4 - f_{24}}{k_2 k_4} = O(\frac{m(N_*)^2
N_{34}}{N_3 (A+N_3)}).$$ But from the mean-value theorem we have
$$ f(k_2) - f(k_{234}) + f(k_3) - f(k_3 - k_{34}) = O(m(N_*)^2 N_{34}),$$
and the claim follows by dividing by $k_2 k_4$.

{\bf Case 4: $N_2 \sim N_3 \sim N_4$.}  Observe this case is essentially symmetric in the indices
$1,2,3,4$. By definition of $M_4$, $\sigma_3$, $\alpha_3$ we have \bas |M_4| &\sim
|[\frac{(f_1 + f_2 + f_{34}) b_{34}}{k_1 k_2}]_{sym}| \\
&\sim N_1^{-4} |[(f_1 + f_2 + f_{34}) b_{34} k_3 k_4]_{sym}|.
\end{align*}
Our task is thus to show that
$$ [(f_1 + f_2 + f_{34}) b_{34} k_3 k_4]_{sym} = O( m(N_*)^2 N_{12} N_{23} N_{13} ).$$
Since $b_{34} = b_{12}$, it will suffice by symmetry to show that
$$ (f_1 + f_2 + f_{34}) k_3 k_4 + (f_3 + f_4 + f_{12}) k_1 k_2 =
O( m(N_*)^2 N_{12} N_{23} N_{13} ).$$ Observe the identity
$$ k_3 k_4 - k_1 k_2 = k_3 k_4 + k_{234} k_2 = k_{23} k_{24}$$
hence we can write the left-hand side as
$$ (f_1 + f_2 + f_3 + f_4) k_1 k_2 + (f_1 + f_2 + f_{34}) k_{23} k_{24}$$
(since $f_{34} = - f_{12}$).  By Lemma \ref{M3-bound}, the second term is $O(m(N_*)^2 N_{34}
N_{23} N_{24})$ which is acceptable.  Thus it will suffice to show that
$$ f_1 + f_2 + f_3 + f_4 = O( m(N_*)^2 N_{12} N_{23} N_{13} / N_1^2).$$
Since $k_{12}+k_{13}+k_{23} = -2k_4$, we see that at least one of $N_{12}, N_{13}, N_{23}$ is
comparable to $N_1$.  Without loss of generality we may take $N_{23} \sim N_1$.  We now write the
left-hand side as
$$ f(k_1) - f(k_1 - k_{12}) - f(k_1 - k_{13}) + f(k_1 - k_{12} - k_{13})$$
and using the double mean-value theorem\footnote{One could object that $f''$ is much larger than
$N_1^{-1}$ near the origin.  However, since we are only evaluating $f$ at points in the annulus
$\{ k: |k| \sim N_1 \}$, we can smooth out $f$ inside this annulus so that $f'' = O(N_1^{-1})$
throughout the interval $\{ k: |k| \lesssim N_1 \}$ without affecting the left-hand side.} (since
$f'' = O(N_1^{-1})$ here), to conclude the argument.
\end{proof}


\begin{thebibliography}{10}


\bibitem{borg:xsb}
J. Bourgain, \emph{Fourier restriction phenomena for certain lattice subsets and applications to
nonlinear evolution equations, Part II}, Geometric and Funct. Anal. \textbf{3} (1993), 209-262.

 \bibitem{borg:squeeze}
J. Bourgain, \emph{Aspects of longtime behaviour of solutions of nonlinear Hamiltonian evolution
equations}, GAFA \textbf{5} (1995), 105--140.

\bibitem{borg:nonsqueeze}
J. Bourgain, \emph{Approximation of solutions of the cubic nonlinear Schr\"odinger equations
 by finite-dimensional equations and nonsqueezing properties},
Int. Math. Res. Notices, \textbf{1994}, no. 2, (1994), 79--90.

\bibitem{borg:measures}
J. Bourgain, \emph{Periodic Korteweg de Vries equation with measures as initial data}, Selecta
Math. (N.S.) \textbf{3} (1997), 115-159.

\bibitem{borg:book}
J. Bourgain, \emph{Global solutions of nonlinear Schr\"odinger equations}, AMS
Publications, 1999.

\bibitem{cct}
M. Christ, J. Colliander, T. Tao, \emph{Illposedness for canonical defocussing equations below the
endpoint regularity}, to appear.

\bibitem{ckstt:2a}
J. Colliander, M. Keel, G. Staffilani, H. Takaoka, T. Tao, \emph{Global well-posedness result for
KdV in Sobolev spaces of negative index}, Elec. J. Diff. Eq. \textbf{2001} (2001) No 26, 1--7.

\bibitem{ckstt:2}
J. Colliander, M. Keel, G. Staffilani, H. Takaoka, T. Tao, \emph{Sharp global well-posedness for
periodic and non-periodic KdV and mKdV on $\R$ and $\T$}, J. Amer. Math. Soc. \textbf{16} (2003),
705--749.

\bibitem{ckstt:gKdV}
J. Colliander, M. Keel, G. Staffilani, H. Takaoka, T. Tao, \emph{Multilinear estimates for
periodic KdV equations, and applications}, Journ. Funct. Analy. \textbf{211} (2004), no. 1,
173--218.

\bibitem{cst}
J. Colliander, G. Staffilani, H. Takaoka, \emph{Global well-posedness of the KdV equation below
$L^2$}, Math Res. Letters \textbf{6} (1999), 755-778.

\bibitem{dickey}
L. Dickey, \emph{Soliton equations and Hamiltonian systems}, World Scientific, 1991.

\bibitem{gardner}
C.S. Gardner, \emph{Korteweg-de Vries equation and generalizations IV}, J. Math. Phys. \textbf{12}
(1971), no. 8, 1548--1551.

\bibitem{gromov}
M. Gromov, \emph{Pseudo-holomorphic curves in symplectic manifolds}, Invent. math., \textbf{82}
(1985), 307--347.


\bibitem{hz}
H. Hofer, E. Zehnder, \emph{A new capacity for symplectic manifolds.}, in Analysis et cetera,
Academic Press (1990), 405--428.  Edited by P. Rabinowitz and E. Zehnder.

\bibitem{hz_book}
H. Hofer, E. Zehnder, \emph{Symplectic Invariants and Hamiltonian Dynamics}, Birkh\"auser Verlag,
1994.

\bibitem{kappeler}
T. Kappeler, P. Topalov, \emph{Global well-posedness of KdV in $H^{-1}(\T, \R)$}, preprint 2003.

\bibitem{kt}
T. Kappeler, P. Topalov, \emph{Global fold structure of the Miura map on $L^2$},
IMRN 2004:39 (2004), 2039--2068.

\bibitem{kpv:kdv}
C. Kenig, G. Ponce, L. Vega, \emph{A bilinear estimate with applications to the KdV equation}, J.
Amer. Math. Soc. \textbf{9} (1996), 573--603.

\bibitem{kpv:counter}
C. Kenig, G. Ponce, L. Vega, \emph{On the ill-posedness of some canonical dispersive equations},
Duke Math. J. \textbf{106} (2001), no 3, 617--633..

\bibitem{kpv:negative}
C. Kenig, G. Ponce, L.Vega, \emph{The Cauchy problem for the Korteweg-de Vries equation in Sobolev
spaces of negative indices}, Duke Math. J., \textbf{71} (1993), 1--21.

\bibitem{kuksin}
S. Kuksin, \emph{Infinite-dimensional symplectic capacities and a squeezing theorem for
Hamiltonian PDE's}, CMP \textbf{167} (1995), 521--552.

\bibitem{kuksin_book}
S. Kuksin, \emph{Analysis of Hamiltonian PDEs}, Oxford Lecture Series in Mathematics and Its
Applications, \textbf{19}, Oxford Univ. Press, 2000.

\bibitem{magri}
F. Magri, \emph{A simple model of the integrable Hamiltonian equation}, J. Math. Phys. 
\textbf{19} (1978), no. 5, 1156--1162.

\bibitem{MiuraTransform}
R. Miura, \emph{Korteweg-de Vries equation and generalizations. I. A remarkable explicit nonlinear
transformation}, J. Mathematical. Phys. \textbf{9} (1968), 1202--1204.


\bibitem{olver}
P. Olver, \emph{Applications of Lie groups to differential equations}, Springer, 1997.

\bibitem{early}
A. Sj\"oberg, \emph{On the Korteweg-de Vries equation:  existence and uniqueness}, J. Math. Anal.
Appl. \textbf{29} (1970), 569--579.

\bibitem{TakaokaTsutsumi}
H. Takaoka, Y. Tsutsumi, \emph{Well-posedness of the Cauchy Problem for the modified KdV equation with periodic boundary condition}, Internat. Math. Res. Notices \textbf{2004}, no. 56,  3009--3040.

\bibitem{tao:xsb}
T. Tao, \emph{Multilinear weighted convolution of $L^2$ functions, and applications to non-linear
dispersive equations}, Amer. J. Math. \textbf{123} (2001), 839--908.

\bibitem{tao:wm1}
T. Tao, \emph{Global regularity of wave maps I. Small critical Sobolev norm in high dimension},
IMRN \textbf{7} (2001), 299--328.

\bibitem{tao:wm2}
T. Tao, \emph{Global regularity of wave maps II. Small energy in two dimensions}, Comm. Math.
Phys. \textbf{224} (2001), 443--544.

\bibitem{faddeev}
V.E. Zakharov, L.D. Faddeev, \emph{The Korteweg-de Vries equation is a completely integrable 
Hamiltonian System}, Funkz. Anal. Priloz. \textbf{5}  (1971), no. 4, 18--27.


\end{thebibliography}
\end{document}